\documentclass{article}
\usepackage{amssymb}
\usepackage{amsmath}

\newtheorem{theorem}{Theorem}

\newtheorem{corollary}[theorem]{Corollary}

\newtheorem{definition}[theorem]{Definition}
\newtheorem{example}[theorem]{Example}

\newtheorem{lemma}[theorem]{Lemma}

\newtheorem{proposition}[theorem]{Proposition}

\newenvironment{proof}[1][Proof]{\noindent\textbf{#1.} }{\ \rule{0.5em}{0.5em}}
\input{tcilatex}

\begin{document}

\author{Maurice A. de Gosson \\
Blekinge Institute of Technology, Karlskrona, SE 37179 (Sweden)\\
and\\
University of Colorado at Boulder, Boulder, CO 80302 (USA)}
\title{THE\ ``SYMPLECTIC\ CAMEL\ PRINCIPLE''\ AND\ SEMICLASSICAL\ MECHANICS\ 
}
\maketitle

\begin{abstract}
We propose a theory of semiclassical mechanics in phase space based on the
notion of quantized symplectic area. The definition of symplectic area makes
use of a deep topological property of \ symplectic mappings, known as the
``principle of the symplectic camel'' which places stringent conditions on
the global geometry of Hamiltonian mechanics. Following this principle,
symplectic mappings --and hence Hamiltonian flows-- are much more rigid than
Liouville's theorem suggests. The dynamical objects of our semiclassical
theory are ``waveforms'', whose definition requires the notion of square
root of de Rham forms. The arguments of these square roots are calculated by
using the properties of a generalized Maslov index. The motion of waveforms
is determined by Hamiltonian mechanics, and the local expressions of these
moving waveforms on configuration space are the usual approximate solutions
of WKB-Maslov theory.
\end{abstract}

\tableofcontents

\section{Introduction}

Non-relativistic physics is governed by two Sciences with distinct domains
of applicability: classical mechanics (\textit{CM}), and quantum mechanics (%
\textit{QM}). The paradigm of \textit{CM}\ is Newton's second law 
\begin{equation*}
m\frac{d^{2}x}{dt^{2}}=F 
\end{equation*}
while that of \textit{QM} is Schr\"{o}dinger's equation 
\begin{equation*}
i\hbar \frac{\partial \Psi }{\partial t}=\widehat{H}\Psi \text{.} 
\end{equation*}
Both equations describe motions (of particles in the first case, and of
waves in the second) in \emph{configuration space} $\mathbb{R}_{x}^{n}$.
However, \textit{CM} and \textit{QM }differ profoundly, both physically and
mathematically. They differ \emph{physically}, because \textit{QM}\emph{\ }%
renounces to the idea of material systems with sharply defined positions and
momenta, and\emph{\ }incorporates instead complex probability amplitudes in
its dynamics. They differ \emph{mathematically}, because while Newton's
second law can immediately be interpreted in terms of phase space variables
using the Hamiltonian formalism, there is no simple and obvious way to
define ``phase space wavefunctions''. On the other hand, one of the most
useful manifestations of \textit{QM}, both in physics and chemistry, is 
\emph{semiclassical mechanics} (\textit{SM}), which applies when the scale
relative to $\hbar $ of certain parameters (e.g. position, time, or mass) in
a system is large. Systems to which\ \textit{SM }applies exhibit a behavior
which is both classical and quantal: while certain quantities (for instance
energy or angular momentum) remain quantized, the \emph{motion} of the
system is governed by \textit{CM}. (\textit{SM} is sometimes described as a
way of doing a simplified path-integral formalism with a focus on
``classical paths''.)

The aim of this article is to present a unifying and mathematically rigorous
theory of \emph{semiclassical mechanics in phase space} based on a deep and
striking topological property of Hamiltonian flows, the non-squeezing
theorem. This theorem --also known as the ``principle of the symplectic
camel''-- says that no Hamiltonian flow will ever be able to squeeze a phase
space ball into a phase space cylinder of smaller radius based on a plane of
conjugate variables $x_{j},p_{j}$. We will use the principle of the
symplectic camel, together with a simple physical postulate related to the
quantization of action, to quantize phase space in such a way that we
recover the usual semiclassical energy levels for integrable systems by a 
\emph{purely topological argument}, without any reference, whatsoever, to
the WKB method or to approximate wavefunctions constructed by other methods.

Our paper consists of two parts, which can be read independently:

\begin{itemize}
\item In the first part (Sections \ref{camelou} and \ref{tou}) we begin by
reviewing the ``principle of the symplectic camel'' (which seems to be
little known by physicists). We define the related notion of symplectic
area, which we then use to quantize energy shells by an appropriate physical
postulate on the \emph{periodic orbits} they carry. This postulate
remarkably leads to the correct ground energy levels for the anisotropic
harmonic oscillator in arbitrary dimensions (it can also be used to derive a
classical form of Heisenberg's inequalities as we have shown in \cite{IOP2}%
). We show that our postulate leads, by a topological argument, to the usual
Keller-Maslov quantization condition 
\begin{equation*}
\frac{1}{2\pi \hbar }\oint_{\gamma }pdx-\frac{1}{4}m(\gamma )\text{ \ is an
integer} 
\end{equation*}
in the integrable case ($\gamma $ a loop on the ``invariant torus''). Our
quantization procedure is actually much more general than those find in the
literature (it quantizes periodic orbits and energy shells), and could thus
be applied with profit to systems exhibiting chaotic behavior.

\item The second part (Sections IV to VII) begins by the study of a simple
example, the one-dimensional harmonic oscillator. It contains in embryonic
form the whole theory which is being further developed in Sections V to VII.
We then proceed to survey the notion of phase of a Lagrangian submanifold,
as defined by Leray \cite{Leray} and the cohomological theory of the Maslov
index which we have developed in \cite{CRAS,JMPA,Wiley}. We are thereafter
able to define our \textquotedblleft waveforms\textquotedblright : they are
phase objects whose phase is expressed in terms of the universal covering of
the Lagrangian submanifold, and whose amplitude is the square root of an
arbitrary \textquotedblleft twisted\textquotedblright\ (or de Rham) form.
Our study of the Maslov index will allow us to assign the proper argument to
these square roots. Our constructions apply whether the underlying manifold
is oriented or not (in contrast with other quantization theories where
orientability is a \textit{sine qua non} requirement, as for instance in 
\cite{Sternin}). Finally, we show that the local expressions on
configuration space of our waveforms, whose motion is Hamiltonian, are just\
the usual WKB wavefunctions. \bigskip 
\end{itemize}

\begin{center}
\textit{Notations and terminology}.
\end{center}

The letter $z$ denotes the generic point $(x,p)$ of the phase space $\mathbb{%
R}^{2n}$ $=\mathbb{R}_{x}^{n}\times \mathbb{R}_{p}^{n}$. We equip $\mathbb{R}%
^{2n}$ with the standard symplectic form $\Omega =d(p\,dx)$: 
\begin{equation*}
\Omega =dp\wedge dx=\sum\limits_{j=1}^{n}dp_{j}\wedge dx_{j}\text{.} 
\end{equation*}
We will denote by $\Lambda (n)$ (resp. $Sp(n)$) the Lagrangian Grassmannian
(resp. the symplectic group) of the symplectic space $(\mathbb{R}%
^{2n},\Omega )$: $\ell \in \Lambda (n)$ if and only if $\ell $ is a $n$%
-dimensional linear subspace of $\mathbb{R}^{2n}$ having the property that $%
\Omega (z,z^{\prime })=0$ for all $z$, $z^{\prime }$. The symplectic group $%
Sp(n)$ consists of all automorphisms $s$ of $\mathbb{R}^{2n}$ such that $%
s^{\ast }\Omega =\Omega $, that is $\Omega (sz,sz^{\prime })=\Omega
(z,z^{\prime })$ for all $z,z^{\prime }$. The universal coverings of $%
\Lambda (n)$ and $Sp(n)$ will be denoted by $\Lambda _{\infty }(n)$ and $%
Sp_{\infty }(n)$.

By definition, a symplectic transformation (or: canonical transformation) is
a diffeomorphism of phase space whose Jacobian matrix belongs to $Sp(n)$ at
every point at which it is defined. Also recall that a Lagrangian
submanifold of phase space $\mathbb{R}_{x}^{n}\times \mathbb{R}_{p}^{n}$ is
a $n$-dimensional submanifold $V\subset \mathbb{R}_{x}^{n}\times \mathbb{R}%
_{p}^{n}$ whose tangent spaces are all Lagrangian planes. Equivalently, $%
\iota _{V}^{\ast }\Omega =0$, $\iota _{V}$ being the inclusion operator $%
V\subset \mathbb{R}_{x}^{n}\times \mathbb{R}_{p}^{n}$.

We will also use in Section \ref{argind} some elementary notations from
singular (co-) chain theory. Let $X$ be a non-empty set, $(G,+)$ an Abelian
group and $p$ an integer $\geq 0$. A $G$-valued $p$-cochain on $X$ is a
mapping $c:X^{p+1}\rightarrow G$. The coboundary of $c$ is the $(p+1)$%
-cochain $\partial c$ defined by 
\begin{equation*}
\partial c(x_{0},...,x_{p+1})=\sum_{j=0}^{p+2}(-1)^{j}c(x_{0},...,\hat{x}%
_{j},...,x_{p+2})
\end{equation*}%
where the cap \textquotedblleft\ \symbol{94}\textquotedblright\ deletes the
term it covers. If $\partial c=0$, $c$ is called a $p$-cocycle; if $%
c=\partial b$ for some $(p-1)$-cochain, it is called a coboundary. We have $%
\partial ^{2}c=0$, hence a coboundary is a cocycle.

\section{Symplectic Camel Quantization\label{camelou}}

Let $B(R)$ be a closed ball in phase space with radius $R$:

\begin{equation*}
B(R)=\left\{ (x,p):|x-x_{0}|^{2}+|p-p_{0}|^{2}\leq R^{2}\right\} 
\end{equation*}
and $Z_{j}(r)$ a cylinder with radius $r$: 
\begin{equation*}
Z_{j}(r)=\left\{ (x,p):(x_{j}-x_{0,j})^{2}+(p_{j}-p_{0,j})^{2}\leq
r^{2}\right\} 
\end{equation*}
($1\leq j\leq n$) based on the $x_{j},p_{j}$ plane (we will call hereafter
the $Z_{j}(r)$ \textit{symplectic cylinders}). Gromov \cite{Gromov} proved
in the mid 1980's that there cannot exist a symplectic transformation
sending $B(R)$ inside $Z_{j}(r)$ unless $R\leq r$. In particular, a
Hamiltonian flow can never squeeze a phase space ball inside a symplectic
cylinder with smaller radius. Gromov's theorem is equivalent statement of
the principle of the symplectic camel :

\begin{proposition}
\label{pro}Let $\Pr_{j}:\mathbb{R}_{x}^{n}\times \mathbb{R}%
_{p}^{n}\longrightarrow \mathbb{R}_{x_{j}}\times \mathbb{R}_{p_{j}}$ be the
projection operator. For every symplectic transformation $f$, we have: 
\begin{equation}
\func{Area}\Pr {}_{j}f(B(R))\geq \pi R^{2}\text{.}  \label{6grobis}
\end{equation}
\end{proposition}

(See the proof in \cite{IOP2} where we used (\ref{6grobis}) to derive a
classical form of Heisenberg's uncertainty relations.)

This result is of course striking, because it seems to contradict the common
conception of Liouville's theorem, which is that under a Hamiltonian flow a
volume in phase space can be made as thin as one likes (cf. Gibbs \cite%
{Gibbs} who calls this the ``principle of extension in phase''; also see the
discussion of Liouville's theorem in Penrose \cite{Penrose}). However, what
is overseen is that the proof of Liouville's theorem only uses the fact that
Hamiltonian flows are divergence free. In fact, Hamiltonian flows consist of
symplectic transformations, and this is a much stronger property than being
just volume preserving as soon as $n>1$. For instance, the statement

\begin{quote}
`` $f$ \textit{is a volume-preserving transformation of phase space}''
\end{quote}

is equivalent to saying that if $(x,p)=f(x^{\prime },p^{\prime })$ then the
Jacobian matrix 
\begin{equation*}
f^{\prime }(z)=\left( 
\begin{array}{cc}
\frac{\partial x}{\partial x^{\prime }} & \frac{\partial x}{\partial
p^{\prime }} \\ 
\frac{\partial p}{\partial x^{\prime }} & \frac{\partial p}{\partial
p^{\prime }}%
\end{array}
\right) 
\end{equation*}
has determinant equal to one, while the statement

\begin{quote}
`` $f$ \textit{is a symplectic transformation of phase space}''
\end{quote}

means that the entries of $f^{\prime }(z)$ satisfy the much more stringent
conditions 
\begin{equation*}
\left\{ 
\begin{array}{c}
\left( \frac{\partial x}{\partial x^{\prime }}\right) ^{T}\frac{\partial p}{%
\partial x^{\prime }}\text{, }\left( \frac{\partial p}{\partial p^{\prime }}%
\right) ^{T}\frac{\partial x}{\partial p^{\prime }}\text{ are
symmetric,\medskip } \\ 
\left( \frac{\partial x}{\partial x^{\prime }}\right) ^{T}\frac{\partial p}{%
\partial p^{\prime }}-\left( \frac{\partial p}{\partial x^{\prime }}\right)
^{T}\frac{\partial x}{\partial p^{\prime }}=I_{n\times n}\text{.}%
\end{array}
\right. 
\end{equation*}

No ``easy'' proofs of Gromov's theorem are known. In Gromov's original paper
and in Hofer-Zehnder \cite{HZ} the reader will find proofs making use of the
theory of pseudo-holomorphic curves. Viterbo gives in \cite{Viterbo} a very
interesting alternative proof using the notion of generating function.

\subsection{Symplectic Area and Periodic Orbits}

Let $\mathcal{D}$ be a subset of $\mathbb{R}_{x}^{n}\times \mathbb{R}%
_{p}^{n} $. We will call \textit{symplectic radius} of $\mathcal{D}$ the
supremum $R_{\max }$ of all $R\geq 0$ such that we can send the phase space
ball $B(R)$ inside $\mathcal{D}$ using a symplectic transformation. We will
call \textit{symplectic area} of $\mathcal{D}$, and denote by $\mathcal{A}(%
\mathcal{D})$ the number $\pi R_{\max }^{2}$: 
\begin{equation*}
\mathcal{A}(\mathcal{D})=\sup_{f\text{ symplectic}}\left\{ \pi
R^{2}:f(B(R))\subset \mathcal{D}\right\} \text{.} 
\end{equation*}
($\mathcal{A}(\mathcal{D})$ is also sometimes called the \textit{symplectic
capacity} of $\mathcal{D}$; see e.g. \cite{HZ}). It immediately follows from
the definition of $\mathcal{A}(\mathcal{D})$ that $A(f(\mathcal{D)})=%
\mathcal{A}(\mathcal{D})$ for every symplectic transformation $f$ of phase
space: symplectic area is thus a \emph{symplectic invariant}.

\textit{Remark}. The notion of symplectic area was first introduced by
Ekeland and Hofer in \cite{EH}; there are other non-equivalent definitions
of symplectic areas/capacities (see e.g. \cite{HZ}).

The principle of the symplectic camel can obviously be restated as 
\begin{equation}
B(R)\subset \mathcal{D\subset }Z_{j}(R)\Longrightarrow \mathcal{A}(\mathcal{D%
})=\pi R^{2}  \label{inclus}
\end{equation}
showing that subsets of phase space with very different shapes and volumes
can have the same symplectic area. Let for instance 
\begin{equation*}
\mathcal{E}:\text{ \ }\sum_{j=1}^{n}\frac{1}{R_{j}^{2}}\left(
p_{j}^{2}+x_{j}^{2}\right) \leq 1\text{.} 
\end{equation*}
be a phase space ellipsoid; we assume that $R_{1}\leq \cdot \cdot \cdot \leq
R_{n}$. (The equation of every ellipsoid in phase space can be put in the
form above by a suitable symplectic change of coordinates.) It follows from
property (\ref{inclus}) that the symplectic area of this ellipsoid is 
\begin{equation}
\mathcal{A}(\mathcal{E)=}\pi R_{1}^{2}\text{.}  \label{capips}
\end{equation}
The symplectic area of a set has --as the terminology suggests-- the
dimension of an area. In the case $n=1$ the symplectic area is in fact just
the usual area: 
\begin{equation}
\mathcal{A}(\mathcal{D})=\left| \int_{\mathcal{D}}dpdx\right| \text{.}
\label{un}
\end{equation}
Notice that the symplectic area of a ball $B(R)$, or of a symplectic
cylinder, is independent of the dimension of the ambient phase space, as it
always is $\pi R^{2}$. The symplectic area and the volume of a ball $B(R)$
in $\mathbb{R}_{x}^{n}\times \mathbb{R}_{p}^{n}$ are related by the formula 
\begin{equation}
\func{Vol}B(R)=\frac{1}{n!}\left[ \mathcal{A}(B(R))\right] ^{n}  \label{A}
\end{equation}
since $B(R)$ has volume $\pi ^{n}R^{2n}/n!$.

\subsection{Symplectic Area and Periodic Orbits}

Consider a bounded domain $\mathcal{D}$ with boundary $\gamma $ in the phase
plane $\mathbb{R}^{n}\times \mathbb{R}^{n}$. Obviously formula (\ref{un})
can be written, using Stoke's theorem as 
\begin{equation*}
\mathcal{A}(\mathcal{D})=\left| \int_{\gamma }pdx\right| 
\end{equation*}
showing that

\begin{center}
``\textit{symplectic area }= \textit{action''}
\end{center}

in the case $n=1$. It turns out --and this is another striking feature of
the principle of the symplectic camel-- that this relation holds in any
dimension. In fact, symplectic area is related to the action of periodic
orbits of Hamiltonian systems. Let us begin by some general considerations.
Consider an infinitely differentiable function (``Hamiltonian'') $H:\mathbb{R%
}_{x}^{n}\times \mathbb{R}_{p}^{n}\longrightarrow \mathbb{R}$ and $%
X_{H}=(\nabla _{p}H,-\nabla _{x}H)$ the associated Hamilton vector field.
The \textit{periodic orbits} of $X_{H}$ are defined as follows: let $(f_{t})$
be the flow of $X_{H}$ and assume that there exists $z$ and $T>0$ such that $%
f_{t}(z)=f_{t+T}(z)$. Then $\gamma (t)=f_{t}(z)$, $0\leq t\leq T$ is a
periodic orbit through $z$. The action of the periodic orbit $\gamma $ is
then the integral 
\begin{equation*}
\oint_{\gamma }pdx=\int_{0}^{T}p(t)dx(t) 
\end{equation*}
where $(x(t),p(t))=f_{t}(z)$.

By definition an ``energy shell'' $\Sigma $ of $H$ is. a non-empty regular
level set of the $H$: 
\begin{equation*}
\Sigma =\left\{ (x,p)\in \mathbb{R}_{x}^{n}\times \mathbb{R}%
_{p}^{n}:H(x,p)=E\right\} \text{.} 
\end{equation*}
Any smooth hypersurface of phase space can of course be viewed as the energy
shell of some Hamiltonian function $H$: it suffices to choose for $H$ any $%
C^{\infty }$ function on $\mathbb{R}_{x}^{n}\times \mathbb{R}_{p}^{n}$
keeping the constant value $E$ in a tubular neighborhood of $\Sigma $. By
definition, a periodic orbit of $\Sigma $ is then a periodic orbit of the
flow determined by $H$, and lying on the energy shell $\Sigma $. Of course,
for this definition to make sense, we have to show that these periodic
orbits are independent of the choice of the Hamiltonian having $\Sigma $ for
energy shell. This follows from the following well-known result:

\begin{lemma}
let $H$ and $K$ be two Hamiltonians, and suppose that there exist two
constants $h$ and $k$ such that 
\begin{equation}
\Sigma =\left\{ z:H(z)=h\right\} =\left\{ z:K(z)=k\right\}  \label{7leven}
\end{equation}
with $\nabla _{z}H\neq 0$ and $\nabla _{z}K\neq 0$ on $\Sigma $. Then the
Hamiltonian vector fields $X_{H}$ and $X_{K}$ have the same periodic orbits
on $\Sigma $.
\end{lemma}

\begin{proof}
It suffices to show that $X_{H}$ and $X_{K}$ have the same integral curves,
up to a reparametrization. Since $\nabla _{z}H(z)\neq 0$ and $\nabla
_{z}K(z)\neq 0$ are both normal to $\Sigma $ at $z$, there exists a function 
$\alpha \neq 0$ on $\Sigma $ such that $X_{K}=\alpha X_{H}$ on $\Sigma $.
Let now $(f_{t})$ and $(g_{t})$ be the flows of $H$ and $K$, respectively,
and define a function $t=t(z,s)$, $s\in \mathbb{R}$ as being the solution of
the ordinary differential problem 
\begin{equation*}
\frac{dt}{ds}=\alpha (f_{t}(z))\text{ \ , \ }t(z,0)=0
\end{equation*}%
(where $z$ is being viewed as a parameter). We claim that 
\begin{equation}
g_{s}(z)=f_{t}(z)\text{ \ \textit{for} \ }z\in \Sigma \text{.}
\label{7claim}
\end{equation}%
In fact, by the chain rule 
\begin{equation*}
\frac{d}{ds}f_{t}(z)=\frac{d}{dt}f_{t}(z)\frac{dt}{ds}=X_{H}(f_{t}(z))\alpha
(f_{t}(z))
\end{equation*}%
that is, since $X_{K}=\alpha X_{H}$: 
\begin{equation*}
\frac{d}{ds}f_{t}(z)=X_{K}(f_{t}(z))
\end{equation*}%
which shows that the mapping $s\longmapsto f_{t(z,s)}(z)$ is a solution of
the differential equation $\dot{z}=X_{H}(z)$ passing through $z$ at time $%
s=t(z,0)=0$. By the uniqueness theorem on solutions of systems of
differential equations, this mapping must be identical to the mapping $%
s\longmapsto g_{s}(z)$; hence (\ref{7claim}). Both Hamiltonians $H$ and $K$
thus have the same periodic orbits.
\end{proof}

The general problem of the existence of periodic orbits on a given energy
shell $\Sigma $ is a very difficult one, which has not yet been completely
solved. We have however the following partial result (see \cite{HZ} and the
references therein):

\begin{proposition}
If the hypersurface $\Sigma $ is the boundary of a compact star-shaped
submanifold of phase space, then it carries at least one periodic orbit.
\end{proposition}

(Recall that a submanifold $M$ of an Euclidean space is called star-shaped
if there exists a point $z\in M$ such that the line segment joining $z$ to
any other point $z^{\prime }\in M$ lies inside $M$.) In particular, the
boundary of every closed convex submanifold thus carries a periodic orbit.

The essential relation between the action of periodic orbits and symplectic
area is given by the following theorem:

\begin{theorem}
\label{7hz}Let $M$ be a compact star-shaped submanifold in phase space.
Then: (1) Every periodic orbit $\gamma $ on $\Sigma =\partial M$ is such
that 
\begin{equation}
\left| \oint_{\gamma }pdx\right| \geq \mathcal{A}(M)  \label{6HZin}
\end{equation}
and: (2) There exists at least one periodic orbit $\gamma _{\min }$ whose
action is the symplectic area of $M$: 
\begin{equation*}
\left| \oint_{\gamma _{\min }}pdx\right| =\mathcal{A}(M)\text{.} 
\end{equation*}
\end{theorem}

(See again Hofer-Zehnder's treatise \cite{HZ} for a proof.)\medskip

\textit{Remark. }We conjecture that the property of the symplectic camel is
the key to a better understanding of not only quantum mechanics, but also of
classical phenomena. Consider, for example, adiabaticity. While it is rather
well understood in one dimension (cf. ``Einstein's pendulum''), one must
take the usual physical statements and ``proofs'' of adiabatic invariance in
higher dimensions with more than a critical eye. The existence of the
symplectic invariant $\mathcal{A}(M)$ shows that symplectic geometry is, in
a sense, a two-dimensional geometry ``projected'' in higher dimensions.
Perhaps a general adiabatic principle could be derived from Proposition \ref%
{pro} by showing that adiabatic invariance in the phase plane is sufficient
for deducing more general results. We also conjecture that the principle of
the symplectic camel might play a fundamental role in thermodynamics and
statistical physics (e.g., Bose-Einstein and Fermi-Dirac statistics).
Viterbo \cite{Viterbo2} has given other interesting physical interpretations
of the principle of the symplectic camel.

\section{Phase Space Quantization\label{tou}}

The property of the symplectic camel discussed above can be used to quantize
phase space in a very simple way. We will, in particular, recover the
correct ground energy level for the $n$-dimensional anisotropic oscillator.

\subsection{A Physical Postulate}

We now make the following postulate of physical nature:\bigskip

\textbf{Minimum symplectic area postulate}: \textit{The only physically
admissible periodic orbits are those which lie on hypersurfaces }$\Sigma $%
\textit{\ which are boundaries of convex subsets }$M$\textit{\ of phase
space with symplectic area at least }$\frac{1}{2}h$\textit{. Moreover, if }$%
A(M)=\frac{1}{2}h$\textit{\ then }$\Sigma $\textit{\ effectively carries a
physically admissible minimal periodic orbit }$\gamma _{\min }$.\bigskip

Notice that in view of Theorem \ref{7hz} such a minimal periodic orbit
satisfies 
\begin{equation}
\oint_{\gamma _{\min }}pdx=\tfrac{1}{2}h  \label{minac}
\end{equation}
hence our postulate implies quantization of action. It actually implies much
more, as we are going to see: because of the principle of the symplectic
camel, it is not a mere restatement of (\ref{minac}). As we remarked in the
Introduction to this paper, our Postulate quantizes energy shells and
periodic orbits, and thus applies beyond integrable systems.

We begin by giving an immediate striking application, by showing that the
minimum symplectic area postulate leads to the correct energy levels of the
anisotropic multi-dimensional harmonic oscillator.

\begin{proposition}
\label{70level}Consider the $n$-dimensional harmonic oscillator with
Hamiltonian: 
\begin{equation}
H=\sum_{j=1}^{n}\frac{1}{2m_{j}}\left( p_{j}^{2}+m_{j}\omega
_{j}^{2}x_{j}^{2}\right)   \label{7cocosc}
\end{equation}%
The minimum symplectic area postulate implies that the ground energy level
of that Hamiltonian is 
\begin{equation}
E_{0}=\sum_{j=1}^{n}\tfrac{1}{2}\hbar \omega _{j}\text{.}  \label{770level}
\end{equation}
\end{proposition}

\textit{Proof}. Let $L$ be the diagonal matrix with diagonal entries $%
(m_{j}\omega _{j})^{-1/2}$. The symplectic change of variables $%
(x,p)\longmapsto (Lx,L^{-1}p)$ changes $H$ into 
\begin{equation*}
H^{\prime }=\sum_{j=1}^{n}\frac{\omega _{j}}{2}(p_{j}^{2}+x_{j}^{2})\text{.} 
\end{equation*}
The change of variables preserving both action integrals and symplectic
areas, it is sufficient to prove the theorem for $H^{\prime }$. Each orbit 
\begin{equation*}
\gamma :\left\{ 
\begin{array}{c}
x_{1}=x_{1}^{\prime }\cos \omega _{1}t+p_{1}^{\prime }\sin \omega _{1}t\text{
, }p_{1}=x_{1}^{\prime }\sin \omega _{1}t-p_{1}^{\prime }\cos \omega _{1}t%
\text{ } \\ 
\cdot \cdot \cdot \cdot \cdot \cdot \cdot \cdot \cdot \cdot \cdot \cdot
\cdot \cdot \cdot \cdot \cdot \cdot \cdot \cdot \cdot \cdot \cdot \cdot
\cdot \cdot \cdot \cdot \cdot \cdot \cdot \cdot \cdot \cdot \\ 
x_{n}=x_{n}^{\prime }\cos \omega _{n}t+p_{n}^{\prime }\sin \omega _{n}t\text{
, }p_{n}=x_{n}^{\prime }\sin \omega _{n}t-p_{n}^{\prime }\cos \omega _{n}t%
\end{array}
\right. 
\end{equation*}
lies, not only on the ellipsoid which is the energy shell of the Hamiltonian 
$H^{\prime }$, but also on each of the symplectic cylinders 
\begin{equation*}
Z_{j}(R_{j})=\left\{ (x,p):x_{j}^{2}+p_{j}^{2}=R_{j}^{2}\right\} 
\end{equation*}
with $R_{j}^{2}=$ $x_{j}^{\prime 2}+p_{j}^{\prime 2}$ and $1\leq j\leq n$.
These cylinders carry periodic orbits, and their symplectic areas must thus
satisfy the conditions 
\begin{equation*}
\mathcal{A}(Z_{j}(R_{j}))=\pi R_{j}^{2}\geq \tfrac{1}{2}h 
\end{equation*}
in view of our postulate. If $\gamma _{\min }$ is a minimal periodic orbit,
it will thus satisfy 
\begin{equation*}
E(\gamma _{\min })=\sum_{j=1}^{n}\tfrac{1}{2}\omega
_{j}R_{j}^{2}=\sum_{j=1}^{n}\tfrac{1}{2}\hbar \omega _{j} 
\end{equation*}
which is the result predicted by standard quantum mechanics.

\subsection{Quantization of Integrable Systems}

Let us next consider a completely integrable system with Hamiltonian $H$.
There are thus $n$ independent constants of the motion $%
F_{1}=H,F_{2},...,F_{n}$ in involution: $\{F_{j},F_{k}\}=0$ for $1\leq
j,k\leq n$. Given an energy shell $\Sigma $ of $H$, through every point $%
z_{0}=(x_{0},p_{0})$ of $\Sigma $ passes a Lagrangian submanifold $V$
carrying the orbits passing through $z_{0}$. Moreover, when $V$ is connected
(which we assume) there exists a symplectic transformation 
\begin{equation}
f:V\longrightarrow (S^{1})^{k}\times \mathbb{R}^{n-k}  \label{7tork}
\end{equation}%
where $(S_{j}^{1})^{k}$ is the product of $k$ unit circles $S_{j}^{1}$, each
lying in some coordinate plane $x_{j},p_{j}$ ($f$ can be constructed using
\textquotedblleft action-angle variables\textquotedblright , see e.g. \cite%
{Arnold,HGoldstein}). The minimum symplectic area postulate imposes a
condition on the energy shells of any Hamiltonian. That condition is that
there should be no periodic orbits with action less than $\tfrac{1}{2}h$,
and that there should exist \textquotedblleft minimal periodic
orbits\textquotedblright\ having precisely $\frac{1}{2}h$ as action. In
fact, we have the following result which ties the minimum symplectic
area/action principle to the Maslov index of loops:

\begin{theorem}
\label{7fundamental}Let $V$ be a Lagrangian submanifold associated to a
Liouville integrable Hamiltonian $H$ and carrying minimal action periodic
orbits. Then we have 
\begin{equation}
\frac{1}{2\pi \hbar }\oint_{\gamma }pdx-\frac{1}{4}m(\gamma )=0
\label{7Maslov}
\end{equation}%
for every loop on $V$.
\end{theorem}

\begin{proof}
Since the actions of loops are symplectic invariants, we can use the
symplectomorphism (\ref{7tork}) to reduce the proof to the case $%
V=(S^{1})^{k}\times \mathbb{R}^{n-k}$. Since the first homotopy group of $V$
is 
\begin{equation*}
\pi _{1}((S^{1})^{k}\times \mathbb{R}^{n-k})\equiv \pi _{1}(S^{1})^{k}\equiv
(\mathbb{Z}^{k},+) 
\end{equation*}
it follows that every loop in $V$ is homotopic to a loop of the type: 
\begin{equation*}
\gamma (t)=(\gamma _{1}(t),...,\gamma _{k}(t),0,...,0)\text{ \ , \ }0\leq
t\leq T 
\end{equation*}
where $\gamma _{j}$ are loops on $S^{1}$: $\gamma _{j}(0)=\gamma _{j}(T)$.
On the other hand, every loop on $S^{1}$ is homotopic to a loop $\varepsilon
_{j}(t)=(\cos \omega _{j}t,\sin \omega _{j}t)$,\ $0\leq t\leq T_{j}$ so that
there must exist positive integers $\mu _{j}$ ($1\leq j\leq n$) such that $%
\mu _{1}T_{1}=\cdot \cdot \cdot =\mu _{k}T_{k}=T$. We can thus identify $%
\gamma _{j}$ with $\mu _{j}\varepsilon _{j}$, the loop $\varepsilon _{j}$
described ``$\mu _{j}$ times'': 
\begin{equation*}
\mu _{j}\varepsilon _{j}(t)=(\cos \omega _{j}t,\sin \omega _{j}t)\text{ \ \ }%
0\leq t\leq T 
\end{equation*}
and it follows that any loop in $V=(S^{1})^{k}\times \mathbb{R}^{n-k}$ is
homotopic to a loop $\gamma =\mu _{1}\varepsilon _{1}+\cdot \cdot \cdot +\mu
_{k}\varepsilon _{k}$. We thus have 
\begin{equation*}
\oint_{\gamma }pdx=\sum_{j=1}^{k}\mu _{j}\oint_{\varepsilon _{j}}p_{j}dx_{j} 
\end{equation*}
and using the same argument as that leading to the proof of formula (\ref%
{770level}) in Proposition \ref{70level}, we must have 
\begin{equation*}
\oint_{\gamma _{j}}p_{j}dx_{j}=\tfrac{1}{2}h\text{ \ \ \ }(1\leq j\leq k) 
\end{equation*}
and hence 
\begin{equation*}
\oint_{\gamma }pdx=\frac{1}{2}\left( \sum_{j=1}^{k}\mu _{j}\right) h\text{.} 
\end{equation*}
Now, the Maslov index of such a loop $\gamma $ in $(S^{1})^{k}\times \mathbb{%
R}^{n-k}$ is by definition 
\begin{equation*}
m(\gamma )=2\sum_{j=1}^{k}\mu _{j} 
\end{equation*}
(see formula (\ref{muj}) in Example \ref{masloopg}, Section \ref{argind})
hence the Keller-Maslov condition (\ref{7Maslov}).
\end{proof}

This result motivates the following definition:

\begin{definition}
A Lagrangian submanifold $V$ is said to be \emph{quantized} if 
\begin{equation}
\frac{1}{2\pi \hbar }\oint_{\gamma }pdx-\frac{1}{4}m(\gamma )\text{ \ is an
integer}  \label{7KM}
\end{equation}%
for every loop $\gamma $ in $V$.
\end{definition}

This definition is of course nothing else than the usual Maslov-Keller
quantization condition \cite{Keller,Maslov,MF}, originating historically
from WKB theory. We arrived to it by purely topological considerations.

\section{Waveforms on the Circle}

We consider in this Section the one dimensional oscillator with Hamiltonian
function 
\begin{equation*}
H=\frac{1}{2}(p^{2}+x^{2})\text{.} 
\end{equation*}
The flow determined by Hamilton's equations for $H$ consisting of the
rotations 
\begin{equation*}
s_{t}=\left( 
\begin{array}{cc}
\cos t & \sin t \\ 
-\sin t & \cos t%
\end{array}
\right) 
\end{equation*}
the phase-space trajectories are thus the circles $S_{r}^{1}=\left\{
|z|=r\right\} $. These circles carry a natural length element denoted, with
the usual abuse of notation, $ds=r\,d\theta $, where $\theta $ is the polar
angle.

\subsection{Position of the Problem}

One wants to define on $S_{r}^{1}$ objects whose vocation is to play the
role of waveforms in phase space, in the sense that their local expressions
are, at best the ``true'' wavefunction, at worst their semiclassical
approximation (i.e. approximations for ``small $\hbar $''). One looks for an
expression of the type 
\begin{equation}
\Psi (z)=e^{\tfrac{i}{\hbar }\varphi (z)}a(z)\sqrt{ds}  \label{(2)}
\end{equation}
where the phase $\varphi $ and the amplitude $a$ are real functions, and $%
\sqrt{ds}$ is supposed to have some well defined meaning . Unfortunately,
one immediately encounters two difficulties when one tries to define $%
\varphi $ and $\sqrt{ds}$. First of all, if one wants the theory to be
consistent with semiclassical mechanics, one must require that the
differential of the phase be the action form: 
\begin{equation}
d\varphi =p\,dx=-r^{2}\sin ^{2}\theta \,d\theta \text{.}  \label{dphase}
\end{equation}
Unfortunately there exists no such function $\varphi $ because the 1-form $%
pdx$ is not exact on $S_{r}^{1}$. We can however define a function $\varphi $
satisfying (\ref{dphase}) on the universal covering $\pi :\mathbb{R}%
\rightarrow S_{r}^{1}$ of $S_{r}^{1}$. That covering is defined by $\pi
(\theta )=r(\cos \theta ,\sin \theta )$ and one immediately checks that 
\begin{equation}
\varphi (\theta )=\frac{r^{2}}{2}(\sin \theta \cos \theta -\theta )
\label{phasedef}
\end{equation}
satisfies (\ref{dphase}). We are thus led to consider $\Psi (z)$ as being an
expression of the type 
\begin{equation*}
\Psi (\theta )=e^{\tfrac{i}{\hbar }\varphi (\theta )}a(\theta )\sqrt{%
r\,d\theta } 
\end{equation*}
where one allows $\theta $ to take any real value, which amounts to define
the candidate for being a phase space wavefunction on the \emph{universal
covering }of the circle. However, there is a second, more serious
obstruction because one does not see how to define unambiguously the square
root $\sqrt{ds}=\sqrt{r\,d\theta }$. The simplest way out of this difficulty
is to decide that one should only consider the (for instance, positive)
square root of the \emph{density} $\left| ds\right| $, that is that we take 
\begin{equation}
\Psi (\theta )=e^{\tfrac{i}{\hbar }\varphi (\theta )}a(\theta )\sqrt{\left|
r\,d\theta \right| }  \label{dens}
\end{equation}
which indeed has a well defined meaning. However, there is a serious rub
with that choice because it leads to the wrong energy levels: since we are
actually interested in a single-valued function on $S_{r}^{1}$, we have to
impose the condition 
\begin{equation}
\Psi (\theta +2\pi )=\Psi (\theta )  \label{singlevalue}
\end{equation}
to the expression (\ref{dens}), which is equivalent to the condition 
\begin{equation*}
\varphi (\theta +2\pi )=\varphi (\theta )-2N\pi \hbar 
\end{equation*}
for some integer $N$. By definition of (\ref{phasedef}) this is in turn
equivalent to $r^{2}=2N\hbar $, which leads to the energy levels $%
E_{N}=N\hbar $, instead of the physically correct $E_{N}=(N+\frac{1}{2}%
)\hbar $.

\subsection{The Need for De Rham Forms}

The way out of these difficulties, and which leads to the correct
quantization conditions, is the inclusion in the theory of de Rham's \cite%
{Rham} ``forms of odd kind'' (also called ``twisted'' or ``pseudo'' -forms
in the literature) are to ordinary differential forms what
``pseudo-vectors'' are to ordinary vectors. By this we mean that the local
expressions of the de Rham forms depend on the orientation of the local
charts used to define them (rigorously speaking, the de Rham forms are just
ordinary forms, but defined on the oriented double cover of the manifold).
In the case of the harmonic oscillator, this leads to the following
constructions. Consider the atlas of $S_{r}^{1}$ consisting of the four
half-circles 
\begin{equation*}
\begin{array}{l}
S_{r,up}^{1}=\left\{ z:\left| z\right| =r,\func{Im}z>0\right\} \text{ ,\ }%
S_{r,down}^{1}=\left\{ z:\left| z\right| =r,\func{Im}z<0\right\} \\ 
S_{r,left}^{1}=\left\{ z:\left| z\right| =r,\func{Re}z<0\right\} \text{ ,\ }%
S_{r,right}^{1}=\left\{ z:\left| z\right| =r,\func{Re}z>0\right\}%
\end{array}%
\end{equation*}
together with the projections $f_{up},f_{down}:(x,p)\longrightarrow x$ and $%
f_{left},f_{right}:(x,p)\longrightarrow p$. The atlas thus defined is \emph{%
not} oriented; for example, the transition function on $S_{r,left}^{1}\cap
S_{r,down}^{1}$ has negative sign. The local expressions of $d\theta $ in
the charts defined above are, respectively 
\begin{equation}
\begin{array}{l}
(d\theta )_{up}=\varepsilon (r^{2}-x^{2})^{-1/2}\,dx\text{ \ , \ }(d\theta
)_{down}=\varepsilon (r^{2}-x^{2})^{-1/2}\,dx \\ 
(d\theta )_{left}=\varepsilon (r^{2}-p^{2})^{-1/2}\,dp\text{ \ , \ }(d\theta
)_{right}=\varepsilon (r^{2}-p^{2})^{-1/2}\,dp%
\end{array}
\label{updown}
\end{equation}
where $\varepsilon =\pm 1$ is the orientation induced from the $x$, $p$ axes
on $S_{up}^{1}$, etc. by the diffeomorphisms $f_{up}$, etc. Thus, if the
axes come equipped with their usual orientations, then $\varepsilon =-1$ for 
$(d\theta )_{up}$ and $(d\theta )_{left}$ and $+1$ for $(d\theta )_{down}$
and $(d\theta )_{right}$, and a change of orientation has the effect of
reversing the sign of $\varepsilon $. The formulas (\ref{updown}), which are
characteristic for de Rham forms, suggest that we \emph{define }the
``argument'' of $d\theta $ by 
\begin{equation}
\arg d\theta =\left\{ 
\begin{array}{l}
m(\theta )\pi \text{ in }S_{r}^{1}\setminus \left\{ \pm r\right\} \text{ }
\\ 
\\ 
(m(\theta )+1)\pi \text{ in }S_{r}^{1}\setminus \left\{ \pm ir\right\}%
\end{array}
\right.  \label{arg}
\end{equation}
where the integer $m(\theta )$ is defined by 
\begin{equation}
m(\theta )=\left[ \theta /\pi \right] +1  \label{Maslov}
\end{equation}
the square brackets meaning ``integer part of''. Notice that a change of
orientation of the frame $x,p$ amounts replacing $m(\theta )$ by $m(\theta
+\pi )=m(\theta )+1$. Formulas (\ref{arg}) and (\ref{Maslov}) allow us to
define the square root of $ds=rd\theta $ in each of the sets $%
S_{r}^{1}\setminus \left\{ \pm r\right\} $ and $S_{r}^{1}\setminus \left\{
\pm ir\right\} $. In fact, 
\begin{equation*}
\left\{ 
\begin{array}{l}
\sqrt{ds}=i^{m(\theta )}\sqrt{\left| r\,d\theta \right| }\text{ \ in \ }%
S_{r}^{1}\setminus \left\{ \pm r\right\} \\ 
\\ 
\sqrt{ds}=i^{m(\theta )+1}\sqrt{\left| r\,d\theta \right| }\text{ \ in \ }%
S_{r}^{1}\setminus \left\{ \pm ir\right\}%
\end{array}
\right. 
\end{equation*}
(Notice that both expressions do not coincide on the overlaps). We are thus
lead to give the following definition of $\Psi (\theta )$: it is the phase
space object whose expression on $S_{r}^{1}\setminus \left\{ \pm r\right\} $
is given by 
\begin{equation}
\Psi _{0}(\theta )=e^{\tfrac{i}{\hbar }\varphi (\theta )}a(\theta
)i^{m(\theta )}\sqrt{\left| r\,d\theta \right| }\text{ }  \label{(3)}
\end{equation}
and on $S_{r}^{1}\setminus \left\{ \pm ir\right\} $ by $\Psi _{1}(\theta
)=i\Psi (\theta )$: 
\begin{equation}
\Psi _{1}(\theta )=e^{\tfrac{i}{\hbar }\varphi (\theta )}a(\theta
)i^{m(\theta )+1}\sqrt{\left| r\,d\theta \right| }  \label{(4)}
\end{equation}
With that definition the single-valuedness condition (\ref{singlevalue})
becomes 
\begin{equation}
\Psi _{j}(\theta +2\pi )=\Psi _{j}(\theta )\text{ \ , \ }j=1,2
\label{single'}
\end{equation}
and is equivalent to $r^{2}=(2N+1)\hbar $, which yields the true energy
levels $E_{N}=(N+\frac{1}{2})\hbar $ predicted by quantum mechanics.

We are going to show that this construction of phase space waveforms can be
extended to \emph{any} physical system to which a Lagrangian submanifold can
be associated. We begin by defining a notion of phase on arbitrary
Lagrangian manifolds generalizing (\ref{phasedef}).

\section{The Lagrangian Phase}

In the rest of this article the letter $V$ will denote a connected
Lagrangian submanifold. Lagrangian manifolds are associated in a natural way
to integrable classical physical systems, and to every quantum system:

\begin{example}
The integrable systems of classical mechanics: $V$ is then topologically an
``invariant torus'', or, more generally a product of $k$ circles and $n-k$
lines.
\end{example}

\begin{example}
Let $\psi (x)=a(x)e^{\frac{i}{\hbar }\Phi (x)}$ where $a$ and $\Phi $ are
defined on some connected open subset of configuration space. The graph $%
V=\left\{ (x,\nabla _{x}\Phi (x)\right\} $ is a Lagrangian submanifold.
\end{example}

\subsection{Definition of the phase}

Consider the universal covering $\pi :$ $\check{V}\longrightarrow V$ of the
Lagrangian submanifold $V$. Since $\check{V}$ is simply connected there
exists a differentiable mapping $\varphi :\check{V}\longrightarrow \mathbb{R}
$ such that 
\begin{equation}
d\varphi (\check{z})=pdx\text{ \ if \ }\pi (\check{z})=(x,p)\text{.}
\label{dephase}
\end{equation}%
We will call, following Leray \cite{Leray}, such a function $\varphi $ a 
\emph{phase} of $V$. The phase can be explicitly constructed in the
following way: choose an \textquotedblleft origin\textquotedblright\ $%
z_{0}\in V$, and identify $\check{z}$ $\in $ $\check{V}$ with the homotopy
classes (with fixed endpoints) of paths in $V$ originating at $z_{0}$; the
projection $\pi (\check{z})$ is then the endpoint $z$ of an element $\gamma
_{z_{0}z}$ of the homotopy class $\check{z}$. A phase function is then given
by the formula 
\begin{equation}
\varphi (\check{z})=\int_{\gamma _{z_{0}z}}p\,dx\text{.}  \label{phase}
\end{equation}%
Clearly the integral only depends on the homotopy class $\check{z}$ of $%
\gamma _{z_{0}z}$ in view of Stoke's theorem, because $\Omega =d(p\,dx)$ is
zero on $V$. Also, 
\begin{equation}
d\varphi (\check{z})=pdx\text{ \ if \ }\pi (\check{z})=(x,p)\text{ .}
\label{difase}
\end{equation}%
We observe that the action of the first homotopy group $\pi _{1}(V)=\pi
_{1}(V,z_{0})$ on $\check{V}$ is reflected by the formula 
\begin{equation}
\varphi (\gamma \check{z})=\varphi (\check{z})+\int_{\gamma }p\,dx\text{ \ }
\label{period}
\end{equation}%
for all $\gamma \in \pi _{1}(V)$. Thus $\varphi $ is defined on $V$ if and
only if all the periods $\int_{\gamma }p\,dx$ of $p\,dx$ vanish, i.e. if $V$
is contractible. We leave it to the Reader to check that formula (\ref{phase}%
) leads to the function (\ref{phasedef}) if we require that $\varphi (0)=0$.

\subsection{The action of Hamiltonian flows on the phase}

Consider a function $H=H(x,p,t)$ defined on some open subset $D\times 
\mathbb{R}_{t}$ of the extended phase space $\mathbb{R}_{x}^{n}\times 
\mathbb{R}_{p}^{n}\times \mathbb{R}_{t}$ . We do not assume here that $H$
has any particular form (for instance ``kinetic energy + potential''), but
only that it is a continuously differentiable function; we also make the
simplifying, but not essential, assumption that the solutions of the
corresponding Hamilton's equations 
\begin{equation*}
\dot{x}=\nabla _{p}H\text{ \ , }\dot{p}=-\nabla _{x}H\text{.} 
\end{equation*}
exist for all times, and are uniquely determined by their values at a time $%
t^{\prime }$. We denote by $\left( f_{t,t^{\prime }}\right) $ the associated
time-dependent flow: $f_{t,t^{\prime }}$ is the symplectic transformation
that takes a point $(x^{\prime },p^{\prime })=(x(t^{\prime }),p(t^{\prime
})) $ to the point $(x,p)=(x(t),p(t))$. When $t^{\prime }=0$, we write
simply $f_{t,0}=f_{t}$. The time-dependent flow satisfies the
Chapman-Kolmogorov relation 
\begin{equation}
f_{t,t^{\prime }}f_{t^{\prime },t^{\prime \prime }}=f_{t,t^{\prime \prime }}%
\text{ }  \label{CP}
\end{equation}
for all times $t$, $t^{\prime }$, $t^{\prime \prime }$.

Suppose that we are given,at some time $t^{\prime }$, a Lagrangian
submanifold $V_{t^{\prime }}$, and select a base point $z_{t^{\prime }}$on $%
V_{t^{\prime }}$. This allows us to define the phase $\varphi (\check{z}%
,t^{\prime })$ of $V_{t^{\prime }}$ by formula (\ref{phase}), with $z_{0}$
replaced by $z_{t^{\prime }}$, $\check{z}$ being an element of the universal
covering $\check{V}_{t^{\prime }}$. The manifold $V_{t}=f_{t,t^{\prime
}}(V_{t^{\prime }})$ is also Lagrangian; defining a base point $z_{t}\in
V_{t}$ by $z_{t}=f_{t,t^{\prime }}(z_{t^{\prime }})$, we identify the
universal covering $\check{V}_{t}$ with $\check{V}_{t^{\prime }}$, defining
the projection $\pi _{t}:\check{V}_{t}\rightarrow V_{t^{\prime }}$ by $\pi
_{t}(\check{z})=z(t)$\ if $\pi _{t^{\prime }}(\check{z})=z(t^{\prime })$.
Denoting by $\varphi (\check{z},t)$ the phase of $V_{t^{\prime }}$, we have:

\begin{proposition}
\label{vite}The phases $\varphi (\check{z},t)$ and $\varphi (\check{z}%
,t^{\prime })$ are related by the formula: 
\begin{equation}
\varphi (\check{z},t)=\varphi (\check{z},t^{\prime })+\int_{z(t^{\prime
})}^{z(t)}p\,dx-H\,dt  \label{PC1}
\end{equation}
where the integral is calculated along the trajectory $s\rightarrow
f_{s,t^{\prime }}(z(t^{\prime }))$ ($t^{\prime }\leq s\leq t$) leading from $%
z(t^{\prime })\in V_{t^{\prime }}$ to $z(t)\in V_{t}$.
\end{proposition}

\begin{proof}
We first note that in view of the relative invariance of the Poincar\'{e}%
-Cartan form (see \cite{Marle}), we have 
\begin{equation}
\varphi (\check{z},t^{\prime })+\int_{z(t^{\prime
})}^{z(t)}p\,dx-H\,dt=\int_{f_{t,t^{\prime }}(\gamma
_{z_{0}z})}p\,dx+\int_{z_{t^{\prime }}}^{z_{t}}pdx-H\,dt  \label{PC2}
\end{equation}
where the integral in the left-hand side is calculated along the trajectory $%
s\rightarrow f_{s,t^{\prime }}(z_{t^{\prime }})$ ($t^{\prime }\leq s\leq t$%
), and $f_{t,t^{\prime }}(\gamma _{z_{0}z})$ is the image in $V_{t}$ by $%
f_{t,t^{\prime }}$ of a path in $V_{t^{\prime }}$ whose homotopy class is $%
\check{z}$. Denoting by $\chi (\check{z},t)$ the left hand-side of (\ref{PC2}%
), we thus have, for fixed $t$: 
\begin{equation*}
d\chi (\check{z},t)=p(t)\,dx(t) 
\end{equation*}
so that $\chi (\check{z},t)$ and $\varphi (\check{z},t)$ differ by a
function $K(t)$ only depending on $t$. Since $\chi (\check{z},t^{\prime
})=\varphi (\check{z},t^{\prime })$, we must have $K=0$.
\end{proof}

\subsection{Phase and generating functions}

The notion of phase of a Lagrangian submanifold is related (as is the action
integral, see \cite{ICP}) to the notion of generating function.

Recall (see for instance Arnol'd\cite{Arnold} or Goldstein \cite{HGoldstein}%
) that a symplectic transformation $f$ is \emph{free} if there exists a
function $W$ defined on twice the configuration space and such that if $%
(x,p)=f(x^{\prime },p^{\prime })$ then 
\begin{equation}
pdx=p^{\prime }dx^{\prime }+dW(x,x^{\prime })\text{.}  \label{libre}
\end{equation}
The function $W$ is then called a \textit{free generating} function (or:
generating function of the second kind) for $f$. When $\ f$ is free, the
relation $(x,p)=f(x^{\prime },p^{\prime })$ uniquely determines $x^{\prime }$
in terms of $x$. In fact, (\ref{libre}) being equivalent to 
\begin{equation}
p=\nabla _{x}W(x,x^{\prime })\text{ \ and \ }p^{\prime }=-\nabla _{x^{\prime
}}W(x,x^{\prime })  \label{libree}
\end{equation}
we have by the implicit function theorem: 
\begin{equation}
\det \frac{\partial (x^{\prime },x)}{\partial (x^{\prime },p^{\prime })}%
=\det \frac{\partial x}{\partial p^{\prime }}\neq 0\text{.}  \label{free}
\end{equation}
Suppose now that $H$ is a Hamiltonian function of the type 
\begin{equation}
H=\sum_{j=1}^{n}\frac{1}{2m_{j}}(p_{j}-A_{j}(x,t))^{2}+U(x,t)\text{.}
\label{maxwell}
\end{equation}
It is then easy to prove (see \cite{ICP}) that there exists $\varepsilon >0$
such that for 
\begin{equation}
0<|t-t^{\prime }|<\varepsilon  \label{epsilon}
\end{equation}
the mappings $f_{t,t^{\prime }}$ are free ($(f_{t,t^{\prime }})$ is the
time-dependent flow determined by $H$). Let now $V_{t^{\prime }}$, $V_{t}$
be as in Proposition \ref{vite}. Keeping initial position and time $%
x^{\prime }$ and $t^{\prime }$ \emph{fixed}, every point $x$ is thus
reached, after time $t-t^{\prime }$, by a \emph{unique} trajectory $\Gamma $
emanating from $x$. Suppose now $z=(x,p)\in V_{t}$. That point is the image
by $f_{t,t^{\prime }}$ of a unique point $z^{\prime }=(x^{\prime },p^{\prime
})\in V_{t^{\prime }}$. The mapping $x\longmapsto x^{\prime }$ thus defined
is a local diffeomorphism, whose inverse we will denote by $f_{t,t^{\prime
}}^{X}$. Thus, by definition, $f_{t,t^{\prime }}^{X}(x^{\prime })$ is the
unique element of $\mathbb{R}_{x}^{n}$ such that 
\begin{equation}
f_{t,t^{\prime }}(x^{\prime },p^{\prime })=(f_{t,t^{\prime }}^{X}(x^{\prime
}),p)\text{ .}  \label{ftx}
\end{equation}
The action integral is then, by definition, the integral of the Poincar\'{e}%
-Cartan form along that trajectory; we note that the function 
\begin{equation}
S(x,x^{\prime };t,t^{\prime })=\int_{x^{\prime },t^{\prime
}}^{x,t}p\,dx-H\,dt\text{ }  \label{generer}
\end{equation}
satisfies Hamilton-Jacobi's equation with initial condition $t=t^{\prime }$: 
\begin{equation}
\frac{\partial S}{\partial t}+H(x,\nabla _{x}S,t)=0\text{ \ , \ }%
S_{x^{\prime }t^{\prime }}(x,x^{\prime };t,t)=0\text{.}  \label{hamac}
\end{equation}

From these considerations we easily get the following consequence of
Proposition \ref{vite}:

\begin{corollary}
Under the assumptions above on the $f_{t,t^{\prime }}$ we have 
\begin{equation}
\varphi (\check{z},t)=\varphi (\check{z},t^{\prime })+S(x,x^{\prime
};t,t^{\prime })  \label{fifri}
\end{equation}
where $\check{z}$ has projection $\pi _{t^{\prime }}(\check{z})=(x^{\prime
},p^{\prime })$ on $V_{t^{\prime }}$, and $(x,p)=f_{t,t^{\prime }}(x^{\prime
},p^{\prime })$. The local expression 
\begin{equation}
\Phi (x,t)=\Phi (x^{\prime },t^{\prime })+S(x,x^{\prime };t,t^{\prime })
\label{locex}
\end{equation}
of $\varphi (\check{z},t)$ satisfies Hamilton-Jacobi's equation.
\end{corollary}

Formula (\ref{fifri}) is an immediate consequence of (\ref{PC1}) and (\ref%
{generer}); formula (\ref{locex}) follows from (\ref{hamac}). (See \cite{ICP}
for a detailed study of the relationship between the action integral and
free generating functions).\medskip

\textit{Remark}. When $H$ is a quadratic homogeneous polynomial in the $%
x_{i},p_{j}$, Euler's identity for homogeneous functions yields 
\begin{equation*}
H=\frac{1}{2}(x\cdot \nabla _{x}H+p\cdot \nabla _{x}p) 
\end{equation*}
hence, using Hamilton's equations: 
\begin{equation}
\varphi (\check{z},t)=\varphi (\check{z},t^{\prime })+\frac{1}{2}(p\cdot
x-p^{\prime }\cdot x^{\prime })-H(z,t)(t-t^{\prime })\text{.}  \label{Euler}
\end{equation}

\section{The Argument\ Index\label{argind}}

The construction of an index generalizing the function $m(\theta )=[\theta
/\pi ]+1$ to arbitrary Lagrangian manifolds is rather technical and will be
done in several steps. We have exposed elsewhere (see \cite%
{IHP,Bullsci,Wiley}) a direct cohomological construction of the argument
index based on previous work by Leray \cite{Leray,Leray2} and the author 
\cite{Wiley,IHP}. We adopt here a more concrete point of view by making use
of Souriau's identification of the Lagrangian Grassmannian with the manifold
of all symmetric unitary matrices (see Souriau's original paper \cite%
{Souriau2} and also \cite{Wiley,ICP}). This approach has the advantage that
it allows straightforward numerical computations and that it does not
require any prior knowledge of chain intersection theory.

\subsection{Maslov and Argument Indices for Paths}

We begin by recalling some results from Lagrangian analysis \cite%
{Wiley,ICP,Leray}.

The ``Souriau mapping'' is the mapping $w:\Lambda (n)\longrightarrow U(n%
\mathbb{)}$ defined by 
\begin{equation}
w(\ell )=u\overline{u}^{\ast }=u(u^{T})\text{ \ if \ }\ell =u(\mathbb{R}%
_{p}^{n})  \label{w}
\end{equation}
where $u\in U(n\mathbb{)}$. This mapping is indeed well defined, because if $%
u(\mathbb{R}_{p}^{n})=u^{\prime }(\mathbb{R}_{p}^{n})$ then $u^{\prime }=uh$
for some $h\in O(n)$ and hence $u^{\prime }\overline{u^{\prime }}^{\ast }=u%
\overline{u}^{\ast }$. The mapping $w$ is in fact a \emph{diffeomorphism},
and hence identifies the Lagrangian Grassmannian $\Lambda (n)$ with the
manifold 
\begin{equation}
W(n\mathbb{)=}\left\{ w\in U(n\mathbb{)},w=w^{T}\right\} \text{.}
\label{W(n)}
\end{equation}
of all symmetric unitary matrices. The universal covering $\Lambda _{\infty
}(n)$ of $\Lambda (n)$ can then identified with the subset 
\begin{equation}
W_{\infty }(n,\mathbb{C})=\left\{ (w,\alpha ):w\in W(n,\mathbb{C)},\det
(w)=e^{i\alpha }\right\}  \label{Wuni(n)}
\end{equation}
of $U(n,\mathbb{C)\times C}$, the covering mapping being the projection $%
(w,\theta )\longmapsto w$. It follows that $\Lambda (n)=W_{\infty }(n)/%
\mathbb{Z}$ and hence $\pi _{1}(\Lambda (n))\equiv (\mathbb{Z},\mathbb{+)}$.
The action of $\pi _{1}(\Lambda (n))$ on $\Lambda _{\infty }(n)\equiv
W_{\infty }(n)$ is given by 
\begin{equation}
\lambda ^{k}\cdot \check{z}=(w,\alpha +2k\pi )\text{.}  \label{k sur W}
\end{equation}
where $\lambda $ is the generator of $\pi _{1}(\Lambda (n))$ whose image in $%
\mathbb{Z}$ is $+1$.

Let us write explicitly these identifications in the case $n=1$. The
manifold $\Lambda (1)$ consists of all straight lines $\ell $ through the
origin in the phase plane $\mathbb{R}^{2}=\mathbb{R}_{x}\mathbb{\times R}%
_{p} $. We will denote by $\ell (\theta )$ the line through the origin whose
angle with the $x$-axis is $\theta +\frac{\pi }{2}$ $(\func{mod}\pi )$: $%
\ell (\theta )$ is thus the direction of the tangent to the unit circle at
the point $e^{i\theta }$. Since $\ell (\theta )=e^{i\theta }(\mathbb{R}_{p})$
the Souriau mapping (\ref{w}) associates to $\ell (\theta )$ the complex
number $w(\theta )=e^{2i\theta }$. It follows that we have the
identifications 
\begin{equation}
\ell (\theta )\equiv e^{2i\theta }\text{ \ and }\ell _{\infty }(\theta
)\equiv (e^{2i\theta },2\theta +2k\pi )\text{, }k\in \mathbb{Z}\text{.}
\label{n=1}
\end{equation}
In particular $\mathbb{R}_{p}$ is identified with $+1$ and $(\mathbb{R}%
_{p})_{\infty }$ with $(1,2k\pi )$.

Consider now the tangent plane $\ell (z)$ at a point $z$ of the Lagrangian
submanifold $V$. The mapping 
\begin{equation}
\ell (.):V\longrightarrow \Lambda (n)\text{ \ , \ }z\mapsto \ell (z)
\label{tangent}
\end{equation}%
is continuous and therefore induces a homomorphism $\ell _{\star }$ from the
first homotopy group of $V$ into that of $\Lambda (n)$. In fact, a base
point $z_{0}\in V$ being chosen once for all, the mapping 
\begin{equation*}
\ell _{\star }:\pi _{1}(V,z_{0})\longrightarrow \pi _{1}(\Lambda (n),\ell
_{0})
\end{equation*}%
(with $\ell _{0}=\ell (z_{0})$) associates to every loop $\gamma $ $:$ $%
[0,1]\rightarrow V$ ($\gamma (0)=\gamma (1)=z_{0}$) the loop $\ell _{\star
}(\gamma )$ of Lagrangian planes defined by $\ell _{\star }(\gamma )(t)=\ell
(\gamma (t))$ , $0\leq t\leq 1$. Using the Souriau identification $\Lambda
(n)\equiv W(n)$ we can associate to the loop $\ell _{\star }(\gamma )$ in $%
\Lambda (n)$ the loop $w_{\ast }\gamma $ in $W(n)$ defined by $w_{\ast
}\gamma (t)=w(\ell _{\star }(\gamma )(t))$. Lifting that loop to $\Lambda
_{\infty }(n)\equiv W_{\infty }(n)$ we get a path 
\begin{equation*}
t\longmapsto (w_{\ast }\gamma (t),\arg \det w_{\ast }\gamma (t))\text{ \ \ ,
\ \ }0\leq t\leq 1
\end{equation*}%
where $\arg \det w_{\ast }\gamma (t)$ is a choice of continuously varying
argument, uniquely determined by a choice of $\arg \det w_{\ast }\gamma (0)$%
. Since $w_{\ast }\gamma (0)=w_{\ast }\gamma (1)$ the quantity 
\begin{equation}
m(\gamma )=\frac{1}{2\pi }\left( \arg \det w_{\ast }\gamma (1)-\arg \det
w_{\ast }\gamma (0)\right)   \label{argmas}
\end{equation}%
must be an integer, only depending on the homotopy class of $\gamma $.
Formula (\ref{argmas}) thus defines a function $m:\pi
_{1}(V,z_{0})\rightarrow \mathbb{Z}$ called \emph{Maslov index for loops}.
The integer $m(\gamma )$ can be intuitively interpreted as follows. Since $%
\pi _{1}(\Lambda (n))\equiv (\mathbb{Z},+)$, $\Lambda (n)$ has a
\textquotedblleft hole\textquotedblright . Now, a loop $\gamma $ in $V$
induces a loop in $\Lambda (n)$, namely the loop $t\mapsto \ell (\gamma
(t))=T_{\gamma (t)}V$, and $m(\gamma )$ is the number of times $\ell _{\ast
}\gamma $ turns around the \textquotedblleft hole\textquotedblright\ in $%
\Lambda (n)$.

\begin{example}
\label{masloopg}Suppose that $V$ is the circle $S^{1}$ in $\mathbb{R}_{x}%
\mathbb{\times R}_{p}$ and $\gamma (t)=e^{2\pi it}$ , $0\leq t\leq 1$. We
have $w_{\ast }\gamma (t)=e^{4\pi it}$, $0\leq t\leq 1$. The argument of $%
w_{\ast }\gamma (t)$ varies from $0$ to $4\pi $ when $t$ goes from $0$ to $1$%
; it follows from definition (\ref{argmas}) that $m(\gamma )=2$. The same
argument shows that if $\gamma =\mu _{1}\varepsilon _{1}+\cdot \cdot \cdot
+\mu _{k}\varepsilon _{k}$ is a loop in $(S^{1})^{k}$, where $\varepsilon
_{j}(t)=e^{2\pi it}$ ($0\leq t\leq 1$) is a loop on the $j$-th circle, then 
\begin{equation}
m(\gamma )=2\sum_{j=1}^{k}\mu _{j}\text{.}  \label{muj}
\end{equation}
\end{example}

The fact that $m(\gamma )$ is an even integer in the example above is not
fortuitous. In fact, Souriau \cite{Souriau} has proved that: 
\begin{equation}
V\text{\textit{\ oriented }}\Longrightarrow m(\gamma )\equiv 0\text{ }\func{%
mod}2\ \text{\textit{for all} }\gamma \in \pi _{1}(V,z_{0})
\label{Maslovorient}
\end{equation}
(see \cite{Wiley} for an algebraic proof of this property, and the
generalization to ``$q$-oriented Lagrangian manifolds''; Dazord \cite{Dazord}
gives a related cohomological definition).

Let us next generalize the notion of Maslov index to arbitrary paths in $%
\Lambda (n)$. Let $\gamma _{z_{0}z}$ be a path in $V$ joining $z_{0}$ to a
point $z$ and $\check{z}$ its homotopy class: $\check{z}$ is an element of
the universal covering $\check{V}$ of $V$. If two paths $\gamma _{z_{0}z}$
and $\gamma _{z_{0}z}^{\prime }$ are homotopic, then so are their images $%
\ell _{\ast }(\gamma _{z_{0}z})$ and $\ell _{\ast }(\gamma _{z_{0}z}^{\prime
})$ in $\Lambda (n)$ by 
\begin{equation*}
\ell (\cdot ):V\ni z\longmapsto T_{z}V\in \Lambda (n)\text{.} 
\end{equation*}
This mapping induces a continuous mapping 
\begin{equation}
\ell _{\infty }(\cdot ):\check{V}\longrightarrow \Lambda _{\infty }(n)
\label{lift}
\end{equation}
which to every $\check{z}\in \check{V}$ with representant $\gamma _{z_{0}z}$
associates the homotopy class $\ell _{\infty }(\check{z})$ of $\ell _{\ast
}(\gamma _{z_{0}z})$; obviously the diagram 
\begin{equation}
\begin{array}{lll}
\check{V} & \overset{\ell _{\infty }(.)}{\longrightarrow } & \Lambda
_{\infty }(n) \\ 
{\tiny \pi }\downarrow &  & \quad \downarrow {\tiny \pi } \\ 
V & \overset{\ell (.)}{\longrightarrow } & \Lambda (n)%
\end{array}
\label{diagram}
\end{equation}
is commutative (the vertical arrows being the covering projections). In view
of the identification $\Lambda (n)\equiv W(n)$ we can associate to $\gamma
_{z_{0}z}$ a unique continuous path $t\mapsto w(t)$ ($t\in \lbrack 0,1]$) in 
$W(n)$ such that $\arg \det w(t)=\alpha (t)$, provided that we have
specified an ``initial argument'' $\alpha (0)$ for $w(0)=\ell (z_{0})$.

We now impose the following rather restrictive condition on the endpoints of
the path $\gamma _{z_{0}z}$: we assume that $z$ is such that 
\begin{equation}
\ell (z_{0})\cap \ell (z)=0  \label{transversale}
\end{equation}
and define an ``argument function'' $m_{0}:\check{V}\longrightarrow \mathbb{R%
}$ by the formula 
\begin{equation}
m_{0}(\check{z})=\frac{1}{2\pi }\left( \alpha (1)-\alpha (0)+i\func{Tr}\func{%
Log}(-w(1)w(0)^{-1}\right) +\frac{n}{2}  \label{Maslov/Souriau}
\end{equation}
where $Tr$ means ``trace of'', and where we define the logarithm by 
\begin{equation}
\func{Log}(-w(1)w(0)^{-1})=\int_{-\infty }^{0}\{\left[ \lambda
I+w(1)(w(0))^{-1}\right] ^{-1}-(\lambda -1)^{-1}I\}\,d\lambda \text{ }
\label{Log}
\end{equation}
($I$ the $n\times n$ identity matrix). The right hand side of (\ref{Log})
makes sense in view of the following characterization of transversality of
Lagrangian planes (see \cite{Wiley,ICP,Leray,Souriau2}):

\begin{lemma}
Let $\ell $ and $\ell ^{\prime }$ be two arbitrary Lagrangian planes, and
set $w=w(\ell )$ and $w^{\prime }=w(\ell ^{\prime })$. The condition $\ell
\cap \mathbb{\ell }^{\prime }=0$ is equivalent to $\det (w(w^{\prime
})^{-1}-I)\neq 0$, that is, to the condition that $w(w^{\prime })^{-1}$ has
no $>0$ eigenvalues.
\end{lemma}

We have:

\begin{proposition}
(1) The function $m_{0}$ is integer-valued. It is locally constant on its
domain of definition $\left\{ z\in V:\ell (z)\cap \mathbb{R}%
_{p}^{n}=0\right\} $; (2) $m_{0}$ coincides with the function defined in (%
\ref{Maslov}) when $V=S_{r}^{1}$ and $z_{0}=+1$; (3) we have for all $\gamma
\in \pi _{1}(V,z_{0})$%
\begin{equation}
m_{0}(\gamma \check{z})=m_{0}(\check{z})+m(\gamma )  \label{jump}
\end{equation}
where $m(\gamma )$ is the Maslov index for loops defined by (\ref{argmas}).
\end{proposition}

\textit{Proof}. (1) We have, by definition of $w$: 
\begin{align*}
\exp (\func{Tr}\func{Log}(-w(1)w(0)^{-1})& =(-1)^{n}\det (w(1)w(0)^{-1}) \\
& =(-1)^{n}\left( \exp (i\alpha (1))-\exp (i\alpha (0)\right) )
\end{align*}
and hence $\exp (2\pi im_{0}(\check{z}))=(-1)^{n}e^{in\pi }=1$ so that $%
m_{0}(\check{z})$ $\in \mathbb{Z}$, as claimed. (2) If $n=1$, $V=S_{r}^{1}$,
and $z_{0}=1$ then $\ell (\theta )\equiv w(1)=e^{2i\theta }$ and $\ \ell
(0)\equiv w(0)=1$. On the other hand the logarithm defined by (\ref{Log}) is
given by in the case $n=1$ by 
\begin{equation}
Log(e^{i\alpha })=i\left( \alpha -2\left[ \tfrac{\alpha +\pi }{\pi }\right]
\pi \right)  \label{log1}
\end{equation}
for $\alpha \neq \pi $ ($\func{mod}2\pi $), hence 
\begin{equation*}
\func{Log}(-w(1)w(0)^{-1})=\func{Log}(-e^{2i\theta })=i\left( 2\theta -2%
\left[ \tfrac{2\theta +2\pi }{2\pi }\right] \pi \right) \text{ } 
\end{equation*}
from which follows that $m_{0}(\check{z})=\left[ \theta /\pi \right] +1$, as
claimed. (3) Let $\check{z}$ be the homotopy class of a path $\gamma
_{z_{0}z}$ and $\gamma $ the homotopy class of a loop $\gamma _{z_{0}z_{0}}$%
. Then $\gamma \check{z}$ is the homotopy class of the concatenation $\gamma
_{z_{0}z_{0}}+\gamma _{z_{0}z}$. Formula (\ref{jump}) follows, by definition
(\ref{Maslov/Souriau}) of the Maslov index for loops.\medskip

The last step in the construction of the complete argument index needs the
properties of the Leray index.

\subsection{The Leray index}

The key to the definition of the Maslov index for paths with endpoints in
general position is the cohomological index defined by Leray \cite%
{Leray,Leray2} in the transversal case, and generalized by the author \cite%
{JMPA} to the non-transversal case. We begin by giving a general definition
of the Leray index. Recall that $\Lambda _{\infty }(n)\equiv W_{\infty }(n)$
is the universal covering of the Lagrangian Grassmannian $\Lambda (n)$.

\begin{definition}
\label{defi}A Leray index on $\left( \Lambda _{\infty }(n)\right) ^{2}$ is a
mapping 
\begin{equation*}
m:\left( \Lambda _{\infty }(n)\right) ^{2}\longrightarrow \mathbb{Z} 
\end{equation*}
having the two following properties: (1) the coboundary of $m$, viewed as a
1-cochain, descends to a $Sp(n)$-invariant cocycle $f$ on $\Lambda (n)$: $%
\partial m=\pi ^{\ast }f$ ($\pi $ the projection $\Lambda _{\infty
}(n)\longrightarrow \Lambda (n)$) (2) $m$ is locally constant on each of the
subsets 
\begin{equation}
\left\{ (\ell _{\infty },\ell _{\infty }^{\prime }):\dim (\ell \cap \ell
^{\prime })=k\right\}  \label{sets}
\end{equation}
\ ($0\leq k\leq n$) of $\left( \Lambda _{\infty }(n)\right) ^{2}$.
\end{definition}

Condition $\partial m=\pi ^{\ast }f$ means that 
\begin{equation}
m(\ell _{\infty },\ell _{\infty }^{\prime })-m(\ell _{\infty },\ell _{\infty
}^{\prime \prime })+m(\ell _{\infty }^{\prime },\ell _{\infty }^{\prime
\prime })=f(\ell ,\ell ^{\prime },\ell ^{\prime \prime })  \label{cobord}
\end{equation}
and the $Sp(n)$-invariance of $f$ means that 
\begin{equation*}
f(s\ell ,s\ell ^{\prime },s\ell ^{\prime \prime })=f(\ell ,\ell ^{\prime
},\ell ^{\prime \prime })\text{ \ for all \ }s\in Sp(n)\text{.} 
\end{equation*}
Notice that the function $f$ automatically is a $\mathbb{Z}$-valued $2$%
-cocycle on $\Lambda (n)$: $\partial f=0$, locally constant on each of the
sets 
\begin{equation}
\left\{ (\ell ,\ell ^{\prime },\ell ^{\prime \prime }):\dim (\ell \cap \ell
^{\prime })=k,\dim (\ell ^{\prime }\cap \ell ^{\prime \prime })=k^{\prime
},\dim (\ell ^{\prime \prime }\cap \ell )=k^{\prime \prime }\right\}
\label{3sets}
\end{equation}
($0\leq k,k^{\prime },k^{\prime \prime }\leq n$). Given a $2$-cocycle $f$ on 
$\Lambda (n)$, there exists at most one Leray index $m$ satisfying (\ref%
{cobord}) (see \cite{JMPA,Wiley}).

We will also need the following simple general property:

\begin{lemma}
\label{lemmuni}Suppose $m$ is a real function defined on all the pairs $%
(\ell _{\infty },\ell _{\infty }^{\prime })$ such that $\ell \cap \ell
^{\prime }=0$, and such that (\ref{cobord}) holds for some $2$-cocycle $f$
on $\Lambda (n)$. Then, the formula 
\begin{equation}
m(\ell _{\infty },\ell _{\infty }^{\prime })=m(\ell _{\infty },\ell _{\infty
}^{\prime \prime })-m(\ell _{\infty }^{\prime },\ell _{\infty }^{\prime
\prime })+f(\ell ,\ell ^{\prime },\ell ^{\prime \prime })  \label{Ledef}
\end{equation}
where $\ell _{\infty }^{\prime \prime }$ is chosen such that $\ell \cap \ell
^{\prime }=\ell \cap \ell ^{\prime \prime }$ defines unambiguously $m(\ell
_{\infty },\ell _{\infty }^{\prime })$ for all $(\ell _{\infty },\ell
_{\infty }^{\prime })\in \left( \Lambda _{\infty }(n)\right) ^{2}$.
\end{lemma}

It is sufficient to verify that $m(\ell _{\infty },\ell _{\infty }^{\prime
}) $ is independent of the choice of $\ell _{\infty }^{\prime \prime }$, but
this is from the cocycle property $\partial f=0$ of $f$ (see \cite%
{JMPA,Wiley,MSDG}).

To every triple $(\ell ,\ell ^{\prime },\ell ^{\prime \prime })$ of
Lagrangian planes we can associate an integer $\sigma (\ell ,\ell ^{\prime
},\ell ^{\prime \prime })$, called \emph{signature}, and defined as being
the difference $\sigma _{+}-\sigma _{-}$ between the number of $>0$ and $<0$
eigenvalues of the quadratic form 
\begin{equation*}
Q(z,z^{\prime },z^{\prime \prime })=\Omega (z,z^{\prime })+\Omega (z^{\prime
},z^{\prime \prime })+\Omega (z^{\prime \prime },z) 
\end{equation*}
on $\ell \oplus \ell ^{\prime }\oplus \ell ^{\prime \prime }$ (see \cite%
{Wiley,ICP,Marle}). The signature is an antisymmetric and $Sp(n)$-invariant
cocycle: $\partial \sigma =0$. Introducing the notation $\dim (\ell ,\ell
^{\prime })=\dim \ell \cap \ell ^{\prime }$ we moreover have 
\begin{equation}
\sigma (\ell ,\ell ^{\prime },\ell ^{\prime \prime })\equiv n+\partial \dim
(\ell ,\ell ^{\prime },\ell ^{\prime \prime })\text{ \ , \ }\func{mod}2\text{%
.}  \label{mod2}
\end{equation}

\begin{theorem}
\label{popo}(1) The function $m$ defined by 
\begin{equation}
m(\ell _{\infty },\ell _{\infty }^{\prime })=\frac{1}{2\pi }\left( \alpha
-\alpha ^{\prime }+i\func{Tr}\func{Log}(-w(w^{\prime })^{-1})\right) +\frac{n%
}{2}  \label{define souriau}
\end{equation}
for $\ell _{\infty }\equiv (w,\alpha )$, $\ell _{\infty }^{\prime }\equiv
(w^{\prime },\alpha ^{\prime })$ with transversal projections: $\ell \cap
\ell ^{\prime }=0$ is the Leray index associated to the cocycle 
\begin{equation}
\func{Inert}(\ell ,\ell ^{\prime },\ell ^{\prime \prime })=\frac{1}{2}\left(
\sigma (\ell ,\ell ^{\prime },\ell ^{\prime \prime })+n+\partial \dim (\ell
,\ell ^{\prime },\ell ^{\prime \prime })\right) \text{.}  \label{Inertia}
\end{equation}
($\func{Inert}$ is called the ``index of inertia'' of $(\ell ,\ell ^{\prime
},\ell ^{\prime \prime })$). (2) That Leray index $m$ has the following
properties: 
\begin{equation}
m(\ell _{\infty },\ell _{\infty }^{\prime })+m(\ell _{\infty }^{\prime
},\ell _{\infty })=n+\dim (\ell ,\ell ^{\prime })\text{ \ , \ }m(\ell
_{\infty },\ell _{\infty })=n\text{ }  \label{partic}
\end{equation}
and the action of $\pi _{1}(\Lambda (n))$ on $m$ satisfies 
\begin{equation}
m(\lambda ^{k}\cdot \ell _{\infty },\lambda ^{k^{\prime }}\cdot \ell
_{\infty }^{\prime })=m(\ell _{\infty },\ell _{\infty }^{\prime
})+k-k^{\prime }  \label{beta}
\end{equation}
where $\lambda $ is the generator of $\pi _{1}(\Lambda (n))$ whose natural
image in $\mathbb{Z}$ is $+1$ (cf. (\ref{k sur W})). (3) For $n=1$ we have,
with the notations (\ref{n=1}): 
\begin{equation}
m(\theta ,\theta ^{\prime })=\left[ \tfrac{\theta -\theta ^{\prime }}{\pi }%
\right] +1\text{.}  \label{mn=1}
\end{equation}
\end{theorem}

\textit{Proof}. We first remark that $\func{Inert}(\ell ,\ell ^{\prime
},\ell ^{\prime \prime })$ always is an integer in view of (\ref{mod2}). We
have shown in \cite{CRAS,JMPA} that the function $\mu $ defined on all $%
\left\{ (\ell _{\infty },\ell _{\infty }^{\prime }):\ell \cap \ell ^{\prime
}=0\right\} $ by $\mu =2m-n$ ($m$ defined by (\ref{define souriau}))
satisfies 
\begin{equation*}
\mu (\ell _{\infty },\ell _{\infty }^{\prime })-\mu (\ell _{\infty },\ell
_{\infty }^{\prime \prime })+\mu (\ell _{\infty }^{\prime },\ell _{\infty
}^{\prime \prime })=\sigma (\ell ,\ell ^{\prime },\ell ^{\prime \prime })%
\text{.} 
\end{equation*}
It follows that 
\begin{equation*}
m(\ell _{\infty },\ell _{\infty }^{\prime })-m(\ell _{\infty },\ell _{\infty
}^{\prime \prime })+m(\ell _{\infty }^{\prime },\ell _{\infty }^{\prime
\prime })=\frac{1}{2}\left( \sigma (\ell ,\ell ^{\prime },\ell ^{\prime
\prime })+n\right) 
\end{equation*}
if the planes $\ell $, $\ell ^{\prime }$, $\ell ^{\prime \prime }$ are
pairwise transverse. Since in this case $\partial \dim (\ell ,\ell ^{\prime
},\ell ^{\prime \prime })=0$, the existence of $m$ follows from Lemma \ref%
{lemmuni}, since $\func{Inert}$ obviously is a $Sp(n)$-invariant cocycle.
Formulas (\ref{Inertia}), (\ref{partic}), (\ref{beta}) are obvious
consequences of (\ref{k sur W}), (\ref{define souriau}) when $\ell \cap \ell
^{\prime }=0$, and of (\ref{cobord}) in the general case since $\ell
_{\infty }$ and $k\cdot \ell _{\infty }$ have same projection $\ell $ on $%
Lag(n)$. Let us finally prove property (3). Suppose first that $\theta
-\theta ^{\prime }\neq 0$ ($\func{mod}$ $\pi $). Then (\ref{mn=1})
immediately follows from (\ref{define souriau}). Suppose next that $\theta
-\theta ^{\prime }=k\pi $. Choosing $\theta ^{\prime \prime }$ such that $%
\ell (\theta )\cap \ell (\theta ^{\prime \prime })=0$, we have 
\begin{equation*}
m(\theta ,\theta ^{\prime })=k+\func{Inert}(\ell (\theta ),\ell (\theta
),\ell (\theta ^{\prime \prime })=k+1 
\end{equation*}
in view of (\ref{Ledef}), concluding the proof.\medskip

\textit{Remark}. There is a deep and interesting connection between the
Leray index $m$ and the Maslov index on the metaplectic group $Mp(n)$ (i.e.
the unitary representation of the double cover $Sp_{2}(n)$ of $Sp(n)$); see 
\cite{AIF,Cocycles}.

\subsection{Definition of the Argument Index}

We now have developed the machinery we need to define the complete argument
index generalizing (\ref{arg})--(\ref{Maslov}).

Consider again the mapping $\ell _{\infty }(.):\check{V}\longrightarrow
\Lambda _{\infty }(n)$ defined by (\ref{diagram}). We denote by $\ell
_{\alpha ,\infty }$ an element of $\Lambda _{\infty }(n)$ with projection $%
\ell _{\alpha }\in \Lambda (n)$.

\begin{proposition}
The function $m_{\alpha }:\check{V}\longrightarrow \mathbb{Z}$ defined by 
\begin{equation}
m_{\alpha }(\check{z})=m(\ell _{\infty }(\check{z}),\ell _{\alpha ,\infty })
\label{mascom}
\end{equation}
has the following properties: (1) Suppose that $\ell _{\alpha }=\ell _{0}$
and $\ell _{\alpha ,\infty }=\ell _{0,\infty }$ is the homotopy class of the
constant loop with origin $\ell _{0}$. Then $m_{\alpha }(\check{z})$ is
given by (\ref{Maslov/Souriau}) if $\ell (z)\cap \ell _{0}=0$; (2) We have 
\begin{equation}
m_{\alpha }(\gamma \check{z})=m_{\alpha }(\gamma \check{z})+m(\gamma )
\label{mapha}
\end{equation}
for all $\gamma \in \pi _{1}(V)$ and $\check{z}\in \check{V}$; (3) We have $%
m_{\alpha }(\check{z})$ $=$ $m(\theta )$ when $V=S^{1}$ , $\ell _{\alpha }=%
\mathbb{R}_{p}$ and $z_{0}=+1$.
\end{proposition}

Property (1) is obvious in view of (\ref{Maslov/Souriau}) and (\ref{define
souriau}). Property (2) follows from property (\ref{beta}) of the Leray
index. Property (3) follows from part (3) of Theorem \ref{popo}.

The following result makes explicit the effect of a change in base point:

\begin{proposition}
Let $m_{\alpha }$, $m_{\beta }$ be the Maslov indices associated by (\ref%
{Maslov/Souriau}) to arbitrary elements $\ell _{\alpha ,\infty }$ and $\ell
_{\beta ,\infty }$ of $\Lambda _{\infty }(n)$. We have 
\begin{equation}
m_{\alpha }(\check{z})-m_{\beta }(\check{z})=m(\ell _{\alpha ,\infty },\ell
_{\beta ,\infty })-\func{Inert}(\ell _{\alpha },\ell _{\beta },\ell (z))
\label{chabase}
\end{equation}
where $z$ is the projection on $V$ of $\check{z}\in \check{V}$.
\end{proposition}

Formula (\ref{chabase}) is of course an immediate consequence of property (%
\ref{cobord}) with the choice $f=\func{Inert}$.\medskip

We next define the waveforms on Lagrangian manifolds.

\section{Waveforms}

We set out to generalize the construction of the square root of a de Rham
form on the circle, as outlined above, to the general case of a $n$%
-dimensional Lagrangian submanifold $V$ (which we again suppose connected,
but \emph{not necessarily orientable}). We will use the following standard
notation and terminology: the \emph{caustic }of $V$ is the closed set \emph{%
\ } 
\begin{equation*}
C=\left\{ z\in V:\ell (z)\cap \mathbb{R}_{p}^{n}\neq 0\right\} \text{ .}
\end{equation*}%
More generally, we will call \textquotedblleft caustic of $V$ relatively to
the direction $\ell _{\alpha }\in \Lambda (n)$\textquotedblright\ the closed
set 
\begin{equation*}
C_{\alpha }=\left\{ z\in V:\ell (z)\cap \ell _{\alpha }\neq 0\right\} 
\end{equation*}%
and we denote by $V_{\alpha }$ its complement $V\setminus C_{\alpha }$: 
\begin{equation*}
V_{\alpha }=\left\{ z\in V:\ell (z)\cap \ell _{\alpha }=0\right\} \text{ .}
\end{equation*}

\subsection{De Rham Forms and their Square Roots}

\textit{Remark}. The introduction of de Rham forms in \textit{SM} should not
be too surprising after all. It is well-known in Physics that many phenomena
exhibit this dependence on orientation, the most elementary example of this
phenomenon being the magnetic field, which is a \textquotedblleft
pseudo-vector\textquotedblright\ (see the lucid discussion in Frankel's book 
\cite{Frankel}). On the other hand, the necessity of the inclusion of 
\textit{half}-densities (or \textit{half}-forms) in quantum mechanics has
been remarked a long time ago (it is of course immediately suggested by Van
Vleck's formula \cite{vanvleck} (also see \cite{ICP,Gutz}). (Historically,
the systematic use of these objects goes back to the work of Blattner,
Kostant and Sternberg (see \cite{Blattner,Kostant,Woodhouse}).)

A \emph{m-density} $\rho \in \left\vert \Omega ^{m}(V)\right\vert $ on $V$ \
($m\in \mathbb{R}$) is a smooth section of the line-bundle $\left\vert
\Lambda ^{m}(TV)\right\vert $ of 1-densities on $TV$. Recall that by
definition every $\rho (z)\in \left\vert \Lambda ^{m}(T_{z}V)\right\vert $
is a mapping 
\begin{equation*}
\rho (z):\underset{m\text{ times}}{\underbrace{T_{z}V\times \cdot \cdot
\cdot \times T_{z}V}}\longrightarrow \mathbb{C}
\end{equation*}%
such that 
\begin{equation*}
\rho (z)(u_{1},...,u_{n})=\left\vert \det A\right\vert ^{m}\rho
(z)(Au_{1},...,Au_{n})
\end{equation*}%
for every $A\in GL(n,\mathbb{C})$ and all vectors $u_{1},...,u_{n}$ in $%
T_{z}V$. Let $(U_{\alpha },f_{\alpha })_{\alpha }$ be an atlas of $V$. The
local expression $\rho _{\alpha }$ of $\rho \in \left\vert \Omega
^{m}(V)\right\vert $ in each chart $(U_{\alpha },f_{\alpha })$ is $\rho
_{\alpha }(x)\left\vert dx\right\vert ^{m}$, the functions $\rho _{\alpha
}\in C^{\infty }(f_{\alpha }(U_{\alpha }))$ satisfy the matching conditions 
\begin{equation}
\rho _{\alpha }(x)=\left\vert \det \frac{\partial f_{\alpha \beta }}{%
\partial x}(x)\right\vert ^{m}\rho _{\beta }(x)\text{ \ , }x\in f_{\beta
}(U_{\alpha }\cap U_{\beta })  \label{match}
\end{equation}%
where we have set $f_{\alpha \beta }=f_{\alpha }\circ f_{\beta }^{-1}$. In
particular, if $f_{\alpha }$ and $f_{\alpha }^{\prime }$ are two coordinate
systems on $U_{\alpha }$, then 
\begin{equation*}
\rho _{\alpha }(x)\left\vert dx\right\vert ^{m}=\left\vert \det \frac{%
\partial x}{\partial x^{\prime }}\right\vert ^{m}\rho _{\alpha }(x^{\prime
})\,\left\vert dx^{\prime }\right\vert ^{m}
\end{equation*}%
if $x=f_{\alpha }(z)$, $x^{\prime }=f_{\alpha }^{\prime }(z)$.

Suppose now that $m=1$; we write $\left| \Omega ^{1}(V)\right| =\left|
\Omega (V)\right| $ and call elements of $\left| \Omega (V)\right| $ simply 
\emph{densities}. Obviously, each $\rho (z)\in \left| \Lambda
^{1}(T_{z}V)\right| $ is homogeneous with respect to scalar multiplication,
but it is not additive. The definition of de Rham forms reinstates
additivity: a \emph{de Rham form} $\tilde{\mu}\in $ $\Omega _{\tau }(V)$
associated to a density $\rho \in \left| \Omega (V)\right| $ is a smooth
section of the line bundle $\tilde{\Lambda}(V)$ obtained by assigning to
each $\rho (z)\in \left| \Lambda (T_{z}V)\right| $ the mapping 
\begin{equation*}
\mu (z):\underset{n\text{ times}}{\underbrace{T_{z}V\times \cdot \cdot \cdot
\times T_{z}V}}\longrightarrow \mathbb{C} 
\end{equation*}
defined by $\mu (z)(u_{1},...,u_{n})=0$ if the vectors $u_{1},...,u_{n}$ are
linearly dependent, and by 
\begin{equation*}
\mu (z)(\mathcal{B}^{\pm }(z))=\pm \rho (z)(\mathcal{B}(z)) 
\end{equation*}
if they form a basis $\mathcal{B}(z)$ of $\ell (z)=T_{z}V$; the notation $%
\pm $ refers to whether this basis is positively or negatively oriented
relatively to the orientation at $z$ defined by a local chart $(U,f)$ at $z$%
. Due to the inclusion of the factor $\pm 1$ in its definition, $\mu (z)$ is
linear, and antisymmetric.

Let us now be more specific, and assume again that $V$ is a Lagrangian
submanifold. We denote by $f_{\alpha }$ the orthogonal projection of $%
V_{\alpha }$ on $\ell _{\alpha }$; it is a local diffeomorphism onto its
image, so that each orientation of $\ell _{\alpha }$ determines an
orientation at $z\in V_{\alpha }$. Let now $U$ be an open neighborhood of $z$
in $V_{\alpha }$. Choosing $U$ sufficiently small, the pair $(U,f_{\alpha })$
is a local chart at $z$. The open set $U$ is orientable, and each of its
orientations is determined by the choice of an orientation of $\ell _{\alpha
}$, that is by the datum of an element $\tilde{\ell}_{\alpha }$of the double
cover $\Lambda _{2}(n)$ with projection $\ell _{\alpha }$. The restriction $%
\ell _{U}(\cdot )$ of the mapping $z\mapsto \ell (z)$ to $U$ can be lifted
to two continuous mappings $z\mapsto \tilde{\ell}_{U}^{\pm }(z)\in \Lambda
_{2}(n)$, corresponding to a continuous positive (resp. negative) choice of
orientations of the tangent planes. Let $\rho \in \left\vert \Omega
\right\vert (U)$ be a density on $U$, and $\ell _{U,\infty }^{\pm }(z)$ be
two elements of $\Lambda _{\infty }(n)$ with projections $\tilde{\ell}%
_{U}^{\pm }(z)\in \Lambda _{2}(n)$. Let $\ell _{\alpha ,\infty }$ be an
element of $\Lambda _{\infty }(n)$ with projection $\tilde{\ell}_{\alpha }$
on $\Lambda _{2}(n)$. We claim that:

\begin{lemma}
The formula 
\begin{equation}
\mu _{U}(z)(\mathcal{B}^{\pm }(z))=(-1)^{m(\ell _{U,\infty }^{\pm }(z),\ell
_{\alpha ,\infty })}\rho (z)(\mathcal{B}(z))  \label{U}
\end{equation}
defines a de Rham form on $U$.
\end{lemma}

In fact, if we change $\mathcal{B}^{+}(z)$ into $\mathcal{B}^{-}(z)$, then
we have to change $\ell _{U,\infty }^{+}(z)$ into $\ell _{U,\infty }^{-}(z)$
in formula (\ref{U}). Since both $\ell _{U,\infty }^{+}(z)$ and $\ell
_{U,\infty }^{-}(z)$ have same projection $\ell (z)\in \Lambda (n)$, we must
have $\ell _{U,\infty }^{+}(z)=\lambda ^{k}\cdot \ell _{U,\infty }^{-}(z)$
for some integer $k$ and hence, by (\ref{beta}) 
\begin{equation*}
m(\ell _{U,\infty }^{+}(z),\ell _{\alpha ,\infty })=m(\ell _{U,\infty
}^{-}(z),\ell _{\alpha ,\infty })+k\text{ .} 
\end{equation*}
Now, the integer $k$ must be odd, because if it where even, then $\ell
_{U,\infty }^{+}(z)$ and $\ell _{U,\infty }^{-}(z)$ would have same
projection $\tilde{\ell}_{U}^{+}(z)$ on $\Lambda _{2}(n)$. Thus 
\begin{equation*}
\mu _{U}(z)(\mathcal{B}^{+}(z))=-\mu _{U}(z)(\mathcal{B}^{-}(z)\text{ .} 
\end{equation*}
Similarly, if we reverse the orientation at $z$, that is, if we replace $%
\tilde{\ell}_{\alpha }$ by an element of $\Lambda _{2}(n)$ defining the
reverse orientation, then we must replace $\ell _{\alpha ,\infty }$ by $%
\lambda \cdot \ell _{\alpha ,\infty }$, which leads to replace $m(\ell
_{U,\infty }^{\pm }(z),\ell _{\alpha ,\infty })$ by 
\begin{equation*}
m(\ell _{U,\infty }^{\pm }(z),\lambda \cdot \ell _{\alpha ,\infty })=m(\ell
_{U,\infty }^{\pm }(z),\ell _{\alpha ,\infty })-1 
\end{equation*}
and thus again reverses the sign of $\mu _{U}(z)(\mathcal{B}^{\pm }(z))$.
The lemma follows, noting that the mappings $z\mapsto \ell _{U,\infty }^{\pm
}(z)$ are locally constant, and hence smooth.

Formula (\ref{U}) allows to define locally the argument of a de Rham form by 
\begin{equation}
\arg \mu _{U}(z)=m(\ell _{U,\infty }^{+}(z),\ell _{\alpha ,\infty }^{+})\pi 
\text{ }  \label{argu}
\end{equation}
and hence the square root of $\mu _{U}$ by the formula 
\begin{equation}
\sqrt{\mu _{U}}(z)(\mathcal{B}^{\pm }(z))=i^{m(\ell _{U,\infty }^{\pm
}(z),\ell _{\alpha ,\infty })}\sqrt{\rho (z)(\mathcal{B}(z))}  \label{locsqr}
\end{equation}

It turns out that this formula can not generally be extended to define a
global argument for a de Rham form (cf. for example the density $r\left|
d\theta \right| $ on the circle $S_{r}^{1}$). However, we are going to show
that this can always be done outside the caustic $C_{\alpha }$ relative to $%
\ell _{\alpha }$, provided that we work on the universal covering of $V$.

\begin{proposition}
Let $\mu $ be a de Rham form on $V$, associated to a density $\rho \in
\left| \Omega \right| (V)$. For every $V_{a}$ there exists a choice of $\ell
_{\alpha ,\infty }\in \Lambda _{\infty }(n)$ such that the restriction $\mu
_{\alpha }$ of $\mu $ to $V_{a}$ is given, for $\check{z}\in $ $\pi
^{-1}(V_{\alpha })$, by 
\begin{equation}
\mu _{\alpha }(\check{z})(\mathcal{B}^{\pm }(z))=(-1)^{m_{\alpha }^{\pm }(%
\check{z})}\rho (z)(\mathcal{B}(z))  \label{rhamalpha}
\end{equation}
where $m_{\alpha }^{+}(\check{z})=m(\ell _{\infty }(\check{z}),\ell _{\alpha
,\infty })$ and $m_{\alpha }^{-}(\check{z})=m(\ell _{\infty }(\check{z}%
),\lambda \cdot \ell _{\alpha ,\infty })$. We can thus define the square
root of $\mu $ of $V_{a}$ by the formula 
\begin{equation}
\sqrt{\mu _{\alpha }}(\check{z})(\mathcal{B}^{\pm }(z))=i^{m_{\alpha }^{\pm
}(\check{z})}\sqrt{\rho (z)(\mathcal{B}(z))}\text{.}  \label{globsqr}
\end{equation}
If $z\in $ $V_{\alpha }\cap V_{\beta }$, then 
\begin{equation}
\sqrt{\mu _{\alpha }}(\check{z})=i^{m_{\alpha \beta }(z)}\sqrt{\mu _{\beta }}%
(\check{z})  \label{transqr}
\end{equation}
where the function $m_{\alpha \beta }:V_{\alpha }\cap V_{\beta }\rightarrow 
\mathbb{Z}$ is given by 
\begin{equation}
m_{\alpha \beta }(z)=m(\ell _{\alpha ,\infty },\ell _{\beta ,\infty })-\func{%
Inert}(\ell _{\alpha },\ell _{\beta },\ell (z))\text{.}  \label{cech}
\end{equation}
\end{proposition}

\textit{Remark}. We have introduced similar notions in \cite{Bullsci,IHP};
however the choice $\mu _{\alpha }=i^{m_{\alpha }(\check{z})}\sqrt{\rho }$
for the square root of a half-density was postulated in a rather ad hoc
manner.

It is instructive to interpret the constructions above in terms of the
oriented double covering $\tilde{V}$ of the manifold $V$. Recall (see for
instance \cite{Godbillon}, X, \S 4) that $\tilde{V}$ is constructed, for an
arbitrary submanifold $V$, in the following way: let $(U_{\alpha },f_{\alpha
})_{\alpha }$ be an atlas, and define, for $U_{\alpha }\cap U_{\beta }\neq
\varnothing $, locally constant mappings $g_{\alpha \beta }:U_{\alpha }\cap
U_{\beta }\longrightarrow \left\{ +1,-1\right\} $ by 
\begin{equation}
g_{\alpha \beta }(z)=Df_{\alpha \beta }(f_{\beta }(z))\left\vert Df_{\alpha
\beta }(f_{\beta }(z))\right\vert ^{-1}\text{.}  \label{transition}
\end{equation}%
where the $f_{\alpha \beta }=f_{\alpha }f_{\beta }^{-1}$ are the transition
functions. Evidently $g_{\alpha \beta }g_{\beta \gamma }g_{\gamma \alpha
}(z)=z$ for $z\in U_{\alpha }\cap U_{\beta }\cap U_{\gamma }$, hence there
exists a twofold covering $\tilde{\pi}:\tilde{V}\longrightarrow V$ with
trivializations 
\begin{equation*}
\tau _{\alpha }:\tilde{\pi}^{-1}(U_{\alpha })\longrightarrow U_{\alpha
}\times \left\{ +1,-1\right\} \text{.}
\end{equation*}%
The transition functions 
\begin{equation*}
\tau _{\alpha \beta }=\tau _{\alpha }\tau _{\beta }^{-1}:U_{\beta }\times
\left\{ +1,-1\right\} \longrightarrow U_{\alpha }\times \left\{
+1,-1\right\} 
\end{equation*}%
are given by $\tau _{\alpha \beta }(z,\varepsilon )=(z,g_{\alpha \beta }(z))$
for $z\in U_{\alpha }\cap U_{\beta }$ and $\varepsilon =\pm 1$. This allows
one to construct an orientable atlas $(\tilde{U}_{\alpha ,\varepsilon },%
\tilde{f}_{\alpha ,\varepsilon })_{\alpha ,\varepsilon }$ of $\tilde{V}$ by
setting $\tilde{U}_{\alpha ,\varepsilon }=\tau _{\alpha }^{-1}(U_{\alpha
}\times \left\{ \varepsilon \right\} )$ and defining $\tilde{f}_{\alpha
,\varepsilon }:\tilde{V}_{\alpha ,\varepsilon }\longrightarrow \mathbb{R}^{n}
$ by the formulas 
\begin{equation*}
\tilde{f}_{\alpha ,\varepsilon }(\tau _{\alpha }^{-1}(z,\varepsilon
))=\left\{ 
\begin{array}{l}
f_{\alpha }(z)\text{ \ for \ }\varepsilon =+1 \\ 
\\ 
\sigma f_{\alpha }(z)\text{ \ for \ }\varepsilon =-1%
\end{array}%
\right. 
\end{equation*}%
where $\sigma $ is $\mathbb{R}^{n}\longrightarrow \mathbb{R}^{n}$ changes,
for instance, the coordinate $x_{1}$ into $-x_{1}$ and leaves the other
coordinates unchanged. One has the following property: $V$\emph{\ is
orientable if and only if the double covering }$\tilde{V}$\emph{\ is
trivial: }$\tilde{V}=V\times \left\{ +1,-1\right\} $\emph{, and }$\tilde{V}$%
\emph{\ is connected if and only }$V$\emph{\ not orientable.}

When $V$ is Lagrangian, we have the following interesting result that shows
that the oriented double cover can always be identified with the product $%
V\times\left\{ +1,-1\right\} $ (but of course not equipped with the product
topology when $V$ is non orientable!):

\begin{proposition}
Suppose that $V$ is not orientable. Then, each of the mappings 
\begin{equation*}
\tilde{\Phi}_{\alpha }:\tilde{V}\longrightarrow V\times \left\{
+1,-1\right\} 
\end{equation*}
defined by $\tilde{\Phi}_{\alpha }(\tilde{z})=(z,(-1)^{m_{\alpha }(\check{z}%
)})$ where $\check{z}\in \check{V}$ has projection $\tilde{z}\in \tilde{V}$,
is a bijection. The restriction $\Phi _{\alpha }$of $\tilde{\Phi}_{\alpha }$
to $\tilde{V}_{\alpha }=\left\{ \tilde{z}:\ell (z)\cap \ell _{\alpha
}\right\} $ is a homeomorphism, and the transitions $\Phi _{\alpha \beta
}=\Phi _{\alpha }\Phi _{\beta }^{-1}$ are given by 
\begin{equation}
\Phi _{\alpha \beta }(\tilde{z})=(z,m(\ell _{\alpha ,\infty },\ell _{\beta
,\infty })-\func{Inert}(\ell _{\alpha },\ell _{\beta },\ell (z))  \label{fab}
\end{equation}
for $\tilde{z}\in \tilde{V}_{\alpha }\cap \tilde{V}_{\beta }$.
\end{proposition}

\begin{proof}
We first note that $\tilde{\Phi}_{\alpha }$ is well-defined: if $\check{z}%
^{\prime }$ and $\check{z}$ both have projection $\tilde{z}\in \tilde{V}$,
then $\check{z}^{\prime }=\gamma \check{z}$ \ for $\gamma \in \pi _{1}(%
\tilde{V})$, and hence 
\begin{equation*}
m_{\alpha }(\check{z}^{\prime })=m_{\alpha }(\gamma \check{z})=m_{\alpha }(%
\check{z})\text{ \ \ }\func{mod}2 
\end{equation*}
in view of (\ref{Maslovorient}) since $\tilde{V}$ is orientable. A similar
argument shows that each mapping $\tilde{\Phi}_{\alpha }$ is injective: if $%
\tilde{\Phi}_{\alpha }(\tilde{z}^{\prime })=\tilde{\Phi}_{\alpha }(\tilde{z}%
) $, then $z^{\prime }=z$ and $\check{z}^{\prime }=\gamma \check{z}$ with $%
m(\gamma )$ even, so that $\gamma \in \pi _{1}(\tilde{V})$, and $\tilde{z}%
^{\prime }=\tilde{z}$. To prove that $\tilde{\Phi}_{\alpha }$ is surjective,
it suffices to note that if $\tilde{\Phi}_{\alpha }(\tilde{z}%
)=(z,\varepsilon )$, then $\tilde{\Phi}_{\alpha }(\tilde{z}^{\prime
})=(z,-\varepsilon )$ where $\tilde{z}^{\prime }$ is the projection on $%
\tilde{V}$ of $\gamma \check{z}$, $\check{z}$ has projection $\tilde{z}$,
and $\gamma \in \pi _{1}(V)$ is such that $m(\gamma )=1$ (the existence of
such a $\gamma $ follows from (\ref{Maslovorient}) since we are assuming $V$
non-orientable). Finally, $\Phi _{\alpha }$ is locally constant on $%
V_{\alpha }$, and hence continuous; formula (\ref{fab}) follows from (\ref%
{chabase}).
\end{proof}

\subsection{Definition of Waveforms}

Let $\varphi $ be the phase of the Lagrangian submanifold $V$, and $\mu $ a
de Rham form associated to a density $\rho $ on $V$.

\begin{definition}
A waveform on $\check{V}$ is the datum, for each $\ell _{\alpha }\in \Lambda
(n)$ of an expression 
\begin{equation*}
\Psi _{\alpha }(\check{z})=e^{\tfrac{i}{\hbar }\varphi (\check{z})}\sqrt{\mu
_{\alpha }}(\check{z})\text{ } 
\end{equation*}
where the $\mu _{\alpha }$ are associated to a same density $\rho $ on $V$.
Equivalently, \ \ 
\begin{equation*}
\Psi _{\alpha }(\check{z})=e^{\tfrac{i}{\hbar }\varphi (\check{z}%
)}i^{m_{\alpha }(\check{z})}\sqrt{\rho }(z) 
\end{equation*}
\end{definition}

Defining the action of $\pi _{1}(V)$ on $\Psi $ by $\gamma \Psi (\check{z}%
)=\Psi (\gamma \check{z})$, we say that $\Psi $ is defined on $V$ if $\gamma
\Psi =\Psi $\ for all \ $\gamma \in \pi _{1}(V)$.

The following result relates our constructions to the minimum symplectic
area postulate:

\begin{proposition}
(1) A waveform is defined on $V$ if and only if $V$ satisfies the Maslov
quantization condition 
\begin{equation}
\dfrac{1}{2\pi \hbar }\int_{\gamma }p\,dx+\frac{1}{4}m(\gamma )\text{ is an
integer}  \label{quantum}
\end{equation}
for every loop $\gamma $ in $V$. (2) When $V$ is oriented, the condition (%
\ref{quantum}) reduces to the condition 
\begin{equation}
\dfrac{1}{\pi \hbar }\int_{\gamma }p\,dx\in \mathbb{Z}\text{ \ for all }%
\gamma \in \pi _{1}(V)\text{.}  \label{orient}
\end{equation}
\end{proposition}

\textit{Proof}. The second statement of the theorem follows from the first
in view of (\ref{Maslovorient}). By definition of a waveform we have to
prove that the condition $\gamma \Psi (\check{z})=\Psi (\gamma \check{z})$
is equivalent to (\ref{quantum}). In view of (\ref{period}) and (\ref{mapha}%
) we have 
\begin{equation*}
\gamma \Psi (\check{z})=\exp \left[ i\left( \dfrac{1}{\hbar }\int_{\gamma
}pdx+\frac{\pi }{2}m(\gamma )\right) \right] \Psi (\check{z}) 
\end{equation*}
hence $\gamma \Psi =\Psi $ is equivalent to 
\begin{equation*}
\dfrac{1}{\hbar }\int_{\gamma }pdx+\frac{\pi }{2}m(\gamma )=0\text{ \ }\func{%
mod}2\pi 
\end{equation*}
which is of course the same thing as (\ref{quantum}).

\subsection{The Hamiltonian motion of waveforms}

The waveforms we have defined are time-independent, they are thus adequate
for the study of stationary processes. However, if we want to use them for
the dynamical study of quantum systems we have to define how Hamiltonian
flows act on them.

Consider an arbitrary function $H=H(x,p,t)$ defined on some open subset $%
D\times\mathbb{R}_{t}$ of the extended phase space $\mathbb{R}_{x}^{n}\times%
\mathbb{R}_{p}^{n}\times\mathbb{R}_{t}$; we denote the time-dependent flow
of $X_{H}=(\nabla_{p}H,-\nabla_{x}H)$ by $f_{t,t^{\prime}}$.

We need the following property of the Leray index. Let $Sp_{\infty }(n)$ be
the universal covering of $Sp(n)$; the usual action $Sp(n)\times \Lambda
(n)\longrightarrow \Lambda (n)$ induces an action $Sp_{\infty }(n)\times
\Lambda _{\infty }(n)\longrightarrow \Lambda _{\infty }(n)$. We claim that
following essential formula holds: 
\begin{equation}
m(s_{\infty }\ell _{\infty },s_{\infty }\ell _{\infty }^{\prime })=m(\ell
_{\infty },\ell _{\infty }^{\prime })  \label{Spin}
\end{equation}
for all $(s_{\infty },\ell _{\infty })\in Sp_{\infty }(n)\times \Lambda
_{\infty }(n)$.

Formula (\ref{Spin}) follows from the fact that the cocycle $f$ associated
to a Leray index is $Sp(n)$-invariant and from the uniqueness \ of the Leray
index associated to a given cocycle (see Lemma \ref{lemmuni}): both mappings 
$(\ell _{\infty },\ell _{\infty }^{\prime })\longmapsto m(\ell _{\infty
},\ell _{\infty }^{\prime })$ and $(\ell _{\infty },\ell _{\infty }^{\prime
})\longmapsto m(s_{\infty }\ell _{\infty },s_{\infty }\ell _{\infty
}^{\prime })$ satisfy the conditions in definition \ref{defi}, and are hence
identical.

The Jacobian matrix $s_{t,t^{\prime }}(z)$ of $f_{t,t^{\prime }}$ at $z$
being symplectic, we can lift the mapping $t\mapsto s_{t,t^{\prime }}(z)\in
Sp(n)$ to a mapping 
\begin{equation*}
t\mapsto (s_{t,t^{\prime }}(z))_{\infty }\in Sp_{\infty }(n) 
\end{equation*}
such that $(s_{t,t}(z))_{\infty }$ is the identity of $Sp_{\infty }(n)$.
Setting 
\begin{equation}
m_{0}(\check{z},t)=m(\ell _{0,\infty },(s_{t,t^{\prime }}(z))_{\infty }\ell (%
\check{z}))  \label{mozart}
\end{equation}
we then define the value of $\Psi $ at time $t$ by the formula: 
\begin{equation}
\Psi (\check{z},t)=e^{\frac{i}{\hbar }\varphi (\check{z},t)}i^{m_{0}(\check{z%
},t)}\sqrt{(f_{t})_{\ast }\rho }(z)\text{.}  \label{timet}
\end{equation}
Let $\check{f}_{t,t^{\prime }}$ be the mapping which to $\Psi (\check{z}%
,t^{\prime })$ associates $\Psi (\check{z},t)$. These mappings satisfy the
Chapman-Kolmogorov condition: 
\begin{equation}
\check{f}_{t,t^{\prime }}\check{f}_{t^{\prime },t^{\prime \prime }}=\check{f}%
_{t,t^{\prime \prime }}\text{ .}  \label{chacha}
\end{equation}
as immediately follows from the fact that 
\begin{equation*}
(f_{t,t^{\prime }}f_{t^{\prime },t^{\prime \prime }})_{\ast }\rho
=(f_{t,t^{\prime }})_{\ast }(f_{t^{\prime },t^{\prime \prime }})_{\ast }\rho 
\text{.} 
\end{equation*}

\subsection{Shadows, and their Motion}

Suppose that the Lagrangian submanifold $V$ is a simply connected graph $%
p=\nabla _{x}\Phi (x)$. $V$ is then automatically quantized, because $%
m(\gamma )=0$ for every loop $\gamma $ in $V$. We denote by $\mathcal{S}$
the operator which to every waveform $\Psi $ on $V$ associates the
coefficient of its local expression in the chart $(V_{t},\pi _{X})$ where $%
\pi _{X}$ is the projection $(x,p)\longmapsto x$ on configuration space.
Thus, if $\Psi $ has local expression $e^{\frac{i}{\hbar }\Phi
(x)}a(x)|dx|^{1/2}$ in $(V,\pi _{X})$, then 
\begin{equation*}
\mathcal{S}(\Psi )(x)=e^{\tfrac{i}{\hbar }\Phi (x)}a(x)\text{ .}
\end{equation*}%
We will call $\Sigma (\Psi )$ the \emph{shadow} of $\Psi $. Suppose now that 
$V$ is an arbitrary Lagrangian submanifold (i.e. we relax the conditions
that $V$ be a graph, or simply connected). We moreover assume that the
quantization condition (\ref{quantum}), holds for $V$. Then, given a point $x
$ there will in general be several charts $(U_{j},\pi _{X})$ such that $x\in
\pi _{X}(U_{j})$. In this case we define the shadow of a waveform $\Psi $ as
being the sum 
\begin{equation*}
\mathcal{S}(\Psi )(x)=\sum_{j}i^{m(x_{j})}e^{\tfrac{i}{\hbar }\Phi
(x,p_{j})}a(x_{j})
\end{equation*}%
calculated at the point $(x,x_{j}^{\prime },t,t^{\prime })$, and $m_{j}$ is
the \emph{Morse index} of the trajectory from $x_{j}^{\prime }$ to $x$: it
is the number of conjugate points along that trajectory, that is, the number
of points where $\det \func{Hess}_{x,x^{\prime }}(S)$ is zero, or infinite.

Write now the wavefunction at time $t^{\prime }$ in the familiar
``oscillatory'' form 
\begin{equation*}
\Psi (x,t^{\prime })=e^{\tfrac{i}{\hbar }\Phi (x,p_{j})}a(x,t^{\prime }) 
\end{equation*}
where $\Phi $ and $a$ are smooth functions, $a\geq 0$, both defined for $%
(x,t^{\prime })\in $ $X\times \mathbb{R}_{t}\ $, $X$ some open subset of $%
\mathbb{R}_{x}^{n}$ (we do not assume that $X$ is simply connected. For
fixed $t^{\prime }$ the function $\psi (\cdot ,t^{\prime })$ is the local
expression of a Lagrangian waveform $\check{\Psi}(.,t^{\prime })$ on the
graph $V_{t^{\prime }}$ of the gradient of the phase $\Phi (.,t^{\prime })$: 
\begin{equation*}
V_{t^{\prime }}=\left\{ (x,p):p=\nabla _{x}\Phi (x,t^{\prime })\right\} 
\text{.} 
\end{equation*}
In fact, 
\begin{equation}
\Psi (\check{z},t^{\prime })=e^{\tfrac{i}{\hbar }\varphi (\check{z}%
,t^{\prime })}\ f^{\ast }\left( a(x,t^{\prime })\left| d^{n}x\right|
^{1/2}\right)  \label{v}
\end{equation}
where $f^{\ast }\left( a(x)\left| d^{n}x\right| ^{1/2}\right) $ is the
pull-back to $V_{t^{\prime }}$ of the half-density $a(x)\left| d^{n}x\right|
^{1/2}$ on $\mathbb{R}_{x}^{n}$.

Let us finally relate our waveforms to the approximate solutions to Schr\"{o}%
dinger's equation studied by Maslov \cite{Maslov} and Maslov and Fedoriuk 
\cite{MF}. Writing the initial wavefunction as $\Psi _{0}(x)=\exp \left[ 
\tfrac{i}{\hbar }S_{0}(x)\right] $, these authors propose expressions of the
type 
\begin{equation}
\Psi (x,t)=\sum_{j}i^{\mu _{j}(x,t,t^{\prime })}\left\vert \frac{dx}{%
dx_{j}^{\prime }}\right\vert ^{-1/2}\exp \left[ \tfrac{i}{\hbar }S_{j}(x,t)%
\right] \Psi _{0}(x_{j}^{\prime },t^{\prime })  \label{magutz}
\end{equation}%
where $x_{j}^{\prime }$, $S_{j}$ and $\mu _{j}(x,t,t^{\prime })$ are defined
as follows: let again $(f_{t,t^{\prime }})$ be the time-dependent flow of $H$%
, and denote by $V_{t}$ the image by $f_{t,t^{\prime }}$ the image of the
Lagrangian submanifold $V_{t^{\prime }}:p=\nabla _{x}S(x,t^{\prime })$.
Given a point $x$ of the projection of $V_{t}$ on $\mathbb{R}_{x}^{n}$ there
is (under adequate assumptions on $U$ and $\Psi _{t^{\prime }}$) a finite
number of points $(x,p_{j})\in V_{t}$ and $(x_{j}^{\prime },p_{j}^{\prime
})\in V_{t^{\prime }}$ such that $(x,p_{j})=f_{t,t^{\prime }}(x_{j}^{\prime
},p_{j}^{\prime })$. The phase $S_{j}(x,t)$ in (\ref{magutz}) is then given
by the integral 
\begin{equation}
S_{j}(x,t)=S(x,x_{j}^{\prime };t,t^{\prime })=\int_{x_{j}^{\prime
},t^{\prime }}^{x,t}pdx-H\,ds  \label{xj}
\end{equation}%
calculated along the trajectory leading from $x_{j}^{\prime }$ at time $%
t^{\prime }$ to $x$ at time $t$. The integers $\mu _{j}(x,t,t^{\prime })$ in
(\ref{magutz}) are the \emph{Morse indices }of these trajectories; these
indices are obtained by counting the number of \emph{conjugate points} along
each trajectory. It turns out that for short time intervals $t-t^{\prime }$
formula (\ref{magutz}) reduces to 
\begin{equation}
\Psi (x,t)=e^{\tfrac{i}{\hbar }S(x,t)}\Psi (x^{\prime },t^{\prime
})\left\vert \det \frac{\partial x}{\partial x^{\prime }}\right\vert ^{-1/2}
\label{a}
\end{equation}%
where $S$ is just the classical action function evaluated from $(x^{\prime
},t^{\prime })$ to $(x,t)$. In fact, if $t-t^{\prime }$ is sufficiently
small then $V_{t}$ will be a graph and there will exist exactly \emph{one}
point $x^{\prime }$ such that $(x,p)=f_{t,t^{\prime }}(x^{\prime },p^{\prime
})$ for $(x,p)\in V_{t}$. (Formula (\ref{a}) was actually already discovered
in 1928 by Van Vleck \cite{vanvleck}.)

\textit{Conclusion}. We have achieved our goal, which was to construct a
semiclassical mechanics based on a topological principle without any
reference to the usual semiclassical approximations based on the WKB method
(from which \textit{SM} historically originates). Semiclassical mechanics
thus appears to be a theory in its own right. We have not examined in which
sense our theory approximates \textit{CM} or \textit{QM}, nor have we given
any applications. As is the case for \textit{CM} or \textit{QM}, the domain
of validity of \textit{SM} can be determined by experience. It would
certainly be interesting to develop further the consequences of the minimum
symplectic area postulate in the following directions:

(1) The quantization of non-integrable systems; it is well-known that
periodic orbits play a fundamental role in such systems (see e.g. \cite%
{BB,Gutz}); the theory of Lagrangian paths as outlined in Section \ref%
{argind} and further developed in \cite{MSDG} is certainly useful in this
context (this theory gives a mathematically rigorous justification of the
recent constructions in Sugita \cite{Sugita2};

(2) Statistical mechanics and thermodynamics: the minimum symplectic area
postulate could be used to push further phase space ``cell quantization'' as
we outlined in \cite{IOP2}. Formula (\ref{A}) relating symplectic area and
volume shows that the volume of such cells corresponding to a quantized ball 
$B(\sqrt{\hbar })$ is $h^{n}/2^{n}n!$, which is consistent with
Bose-Einstein statistics.

\end{document}